\documentclass{cambridge6Atight}
\usepackage{latexsym,amsfonts,amsmath,amssymb,verbatim,txfonts}
\usepackage[numbers]{natbib}
\usepackage{hyperref}
\usepackage{url}

\newcommand{\V}[3]{{}^{\phantom{(.)}}_{#1\!}V_{#2}^{(#3)}}

\newcommand{\R}{\mathbb R}
\newcommand{\abs}[1]{\lvert#1\rvert}

\newcommand{\dup}{\textup{d}}
\newcommand{\eup}{\textup{e}}
\newcommand{\iup}{\textup{i}}

\renewcommand{\theequation}{\thesection.\arabic{equation}}

\numberwithin{equation}{section}

\allowdisplaybreaks

\begin{document}
\pagenumbering{arabic}
\setcounter{page}{1}
\chapterauthor{Hjalmar Rosengren and S.\ Ole Warnaar}
\chapter*{Elliptic hypergeometric functions associated with root systems}

\setcounter{chapter}{1}
\section{Introduction}

Let $f=\sum_{n\geq 0} c_n$.
The series $f$ is called hypergeometric if the ratio $c_{n+1}/c_n$,
viewed as a function of $n$, is rational.
A simple example is the Taylor series $\exp(z)=\sum_{n=0}^{\infty} z^n/n!$.
Similarly, if the ratio of consecutive terms of $f$ is a rational function
of $q^n$ for some fixed $q$ --- known as the base --- then $f$ is called a
basic hypergeometric series.
An early example of a basic hypergeometric series is Euler's $q$-exponential 
function $e_q(z)=\sum_{n\geq 0} z^n/\big((1-q)\cdots(1-q^n)\big)$.
If we express the base as $q=\exp(2\pi\iup/\omega)$ then $c_{n+1}/c_n$
becomes a trigonometric function in $n$, with period $\omega$. 
This motivates the more general definition of an elliptic hypergeometric
series as a series $f$ for which $c_{n+1}/c_n$ is a doubly-periodic
meromorphic function of $n$.

Elliptic hypergeometric series first appeared in 1988 in the work of Date 
et al.\ on exactly solvable lattice models in statistical 
mechanics~\cite{date88}.
They were formally defined and identified as mathematical objects of 
interest in their own right by Frenkel and Turaev in 1997~\cite{ft97}.
Subsequently, Spiridonov introduced the elliptic beta integral, initiating
a parallel theory of elliptic hypergeometric integrals \cite{spir01}.
Together with Zhedanov \cite{spir04,sz00} he also showed that Rahman's 
\cite{rahman86} and Wilson's \cite{wilson91} theory of biorthogonal rational 
functions --- itself a generalization of the Askey scheme \cite{koekoek10} 
of classical orthogonal polynomials --- can be lifted to the elliptic level.

All three aspects of the theory of elliptic hypergeometric functions 
(series, integrals and biorthogonal functions) have been generalized to
higher dimensions, connecting them to root systems and 
Macdonald--Koornwinder theory.
In \cite{warnaar02} Warnaar introduced elliptic hypergeometric series
associated to root systems, including a conjectural series evaluation of
type $\mathrm{C}_n$. 
This was recognized by van Diejen and Spiridonov
\cite{vandspir00,vandspir01} as a discrete analogue of a multiple elliptic
beta integral (or elliptic Selberg integral).
They formulated the corresponding integral evaluation, again as a conjecture.
This in turn led Rains \cite{rains06,rains10} to develop an elliptic
analogue of Macdonald--Koornwinder theory, resulting in continuous as well
as discrete biorthogonal elliptic functions attached to the non-reduced root 
system $\mathrm{BC}_n$. 
In this theory, the elliptic multiple beta integral and its discrete 
analogue give the total mass of the biorthogonality measure.

Although a relatively young field, the theory of elliptic hypergeometric 
functions has already seen some remarkable applications.
Many of these involve the multivariable theory.
In 2009, Dolan and Osborn showed that supersymmetric indices of 
four-dimensional supersymmetric quantum field theories are expressible 
in terms of elliptic hypergeometric integrals \cite{dolosb09}.
 Conjecturally, such field theories admit electric--magnetic dualities
known as Seiberg dualities, such that dual theories have the same index.
This leads to non-trivial identities between elliptic hypergeometric 
integrals (or, for so called confining theories, to integral evaluations). 
In some cases these are known identities, which thus gives a partial 
confirmation of the underlying Seiberg duality.
However, in many cases it leads to new identities that are yet to be
rigorously proved, see e.g.\
\cite{gadde10,gadde10b,spirvart11,spirvart12,spirvart14} and the recent
survey \cite{rastelli17}.
Another application of elliptic hypergeometric functions is to exactly
solvable lattice models in statistical mechanics. 
We already mentioned the occurrence of elliptic hypergeometric series
in the work of Date et al., but more recently it was shown that elliptic
hypergeometric integrals are related to solvable lattice models with
continuous spin parameters \cite{bazser12a,spir11}.
In the one-variable case, this leads to a generalization of many
well-known discrete models such as the two-dimensional Ising model
and the chiral Potts model.
This relation to solvable lattice models has been extended to multivariable
elliptic hypergeometric integrals in \cite{bks13,bazser12b,spir11}.
Further applications of multivariable elliptic hypergeometric functions
pertain to elliptic Calogero--Sutherland-type systems
\cite{razamat14,spir07b} and the representation theory of elliptic
quantum groups \cite{ros10}.

\medskip

In the current chapter we give a survey of elliptic hypergeometric
functions associated with root systems, comprised of three main parts.
The first two form in essence an annotated table of the main evaluation
and transformation formulas for elliptic hypergeometric integrals and
series on root systems. 
The third and final part gives an introduction to Rains' elliptic
Macdonald--Koornwinder theory (in part also developed by
Coskun and Gustafson \cite{coskgust06}).
Due to space limitations, applications will not be covered here and we
refer the interested reader to the above-mentioned papers and references
therein.

Rather than throughout the text, references for the main results 
are given in the form of separate notes at the end of each section.
These notes also contain some brief historical comments and further pointers 
to the literature.

\medskip

\noindent
\textbf{Acknowledgements:} We thank Ilmar Gahramanov, 
Eric Rains, Michael Schloss\-er and Vyacheslav Spiridonov for 
valuable comments.

\subsection{Preliminaries}

Elliptic functions are doubly-periodic meromorphic functions on $\mathbb{C}$.
That is, a meromorphic function $g:\mathbb{C}\to\mathbb{C}$ is elliptic if
there exist $\omega_1,\omega_2$ with $\mathrm{Im}(\omega_1/\omega_2)>0$ 
such that $g(z+\omega_1)=g(z+\omega_2)=g(z)$ for all $z\in\mathbb{C}$.
If we define the elliptic nome $p$ by $p=\eup^{2\pi\iup\omega_1/\omega_2}$
(so that $\abs{p}<1$) then $z\mapsto \eup^{2\pi\iup z/\omega_2}$ maps the
period parallelogram spanned by $\omega_1,\omega_2$ to an annulus with
radii $\abs{p}$ and $1$.
Given an elliptic function $g$ with periods $\omega_1$ and $\omega_2$,
the function $f:\mathbb{C}^{\ast}\to\mathbb{C}$ defined by
\[
g(z)=f\big(\eup^{2\pi\iup z/\omega_2}\big)
\]
is thus periodic in an annulus:
\[
f(pz)=f(z).
\]
By mild abuse of terminology we will also refer to such $f$ as an elliptic
function.
A more precise description would be elliptic function in multiplicative form.

The basic building blocks for elliptic hypergeometric functions are
\begin{align*}
\theta(z)=\theta(z;p)&=\prod_{i=0}^{\infty}(1-zp^i)(1-p^{i+1}/z), \\
(z)_k=(z;q,p)_k&=\prod_{i=0}^{k-1}\theta(zq^i;p), \\
\Gamma(z)=\Gamma(z;p,q)&=
\prod_{i,j=0}^{\infty}\frac{1-p^{i+1}q^{j+1}/z}{1-zp^iq^j}, 
\end{align*}
known as the modified theta function, elliptic shifted factorial and 
elliptic gamma function, respectively. Note that the dependence on the
elliptic nome $p$ and base $q$ will mostly be suppressed from our notation.
One exception is the $q$-shifted factorial 
$(z;q)_{\infty}=\prod_{i\geq 0}(1-zq^i)$ which, to avoid possible confusion,
will never be shortened to $(z)_{\infty}$.

For simple relations satisfied by the above three functions we refer 
the reader to \cite{gr04}. Here we only note that the elliptic gamma
function is symmetric in $p$ and $q$ and satisfies
\[
\Gamma(pq/z)\Gamma(z)=1 \quad \text{and} \quad
\Gamma(qz)=\theta(z)\Gamma(z).
\]

For each of the functions $\theta(z)$, $(z)_k$ and $\Gamma(z)$, we 
employ condensed notation as exemplified by
\begin{align*}
\theta(z_1,\dots,z_m)&=\theta(z_1)\dotsm\theta(z_m), \\
\big(az^{\pm}\big)_k&=(az)_k(a/z)_k,\\[0.4mm]
\Gamma\big(tz^{\pm}w^{\pm}\big)&=\Gamma(tzw)\Gamma(tz/w)
\Gamma(tw/z)\Gamma(t/zw).
\end{align*}

In the trigonometric case $p=0$ we have $\theta(z)=1-z$, so that 
$(z)_k$ becomes a standard $q$-shifted factorial and $\Gamma(z)$ a 
rescaled version of the $q$-gamma function.

We also need elliptic shifted factorials indexed by partitions.
A partition $\lambda=(\lambda_1,\lambda_2,\dots)$ is a weakly decreasing 
sequence of non-negative integers such that only finitely many $\lambda_i$
are non-zero. The number of positive $\lambda_i$ is called the length
of $\lambda$ and denoted by $l(\lambda)$.
The sum of the $\lambda_i$ will be denoted by $\abs{\lambda}$. 
The  diagram of $\lambda$ consists of the points $(i,j)\in\mathbb{Z}^2$
such that $1\leq i\leq l(\lambda)$ and $1\leq j\leq \lambda_i$. 
If these inequalities hold for $(i,j)\in\mathbb{Z}^2$ we write 
$(i,j)\in\lambda$.
Reflecting the  diagram in the main diagonal yields the conjugate
partition $\lambda'$.
In other words, the rows of $\lambda$ are the columns of $\lambda'$ and
vice versa.
A standard statistic on partitions is
\[ 
n(\lambda)=\sum_{i\geq 1} (i-1)\lambda_i
=\sum_{i\geq 1} \binom{\lambda'_i}{2}.
\]
For a pair of partitions $\lambda,\mu$ we write $\mu\subset\lambda$ if
$\mu_i\leq\lambda_i$ for all $i\geq 1$.
In particular, when $l(\lambda)\leq n$ and $\lambda_i\leq N$ for all
$1\leq i\leq N$ we write $\lambda\subset (N^n)$.
Similarly, we write $\mu\prec\lambda$ if the interlacing conditions
$\lambda_1\geq\mu_1\geq\lambda_2\geq\mu_2\geq\cdots$ hold.

With $t$ an additional fixed parameter, we will need the following three 
types of elliptic shifted factorials index by partitions:
\begin{align*}
(z)_{\lambda}=(z;q,t;p)_{\lambda}&=\prod_{(i,j)\in\lambda}
\theta\big(zq^{j-1}t^{1-i}\big)
=\prod_{i\geq 1}\big(zt^{1-i}\big)_{\lambda_i}, \\
C^{-}_{\lambda}(z)=C^{-}_{\lambda}(z;q,t;p)&=\prod_{(i,j)\in\lambda}
\theta\big(zq^{\lambda_i-j}t^{\lambda'_j-i}\big), \\
C^{+}_{\lambda}(z)=C^{+}_{\lambda}(z;q,t;p)&=\prod_{(i,j)\in\lambda}
\theta\big(zq^{\lambda_i+j-1}t^{2-\lambda'_j-i}\big).
\end{align*}
By $\theta(pz)=-z^{-1} \theta(z)$ it follows that $(a)_{\lambda}$
is quasi-periodic:
\begin{equation}\label{quasiperiodicity}
\big(p^k z\big)_{\lambda}=\Big[(-z)^{-\abs{\lambda}} 
q^{-n(\lambda')} t^{n(\lambda)}\Big]^k 
p^{-\binom{k}{2}\abs{\lambda}}\, (z)_{\lambda},
\qquad k\in\mathbb{Z}.
\end{equation}
Again we use condensed notation so that, for example,
$(a_1,\dots,a_k)_{\lambda}=(a_1)_{\lambda}\cdots(a_k)_{\lambda}$. 

\subsection{Elliptic Weyl denominators}

Suppressing their $p$-dependence we define
\begin{align*}
{\Delta^{\mathrm{A}}(x_1,\dots,x_{n+1})}
&={\prod_{1\leq i<j\leq n+1}x_j\, \theta(x_i/x_j),} \\
\Delta^{\mathrm{C}}(x_1,\dots,x_n)&=\prod_{j=1}^n\theta(x_j^2)
\prod_{1\leq i<j\leq n}x_j\,\theta(x_ix_j^{\pm}),
\end{align*}
which are essentially the Weyl denominators of the affine root systems
$\mathrm{A}_n^{(1)}$ and $\mathrm{C}_n^{(1)}$ \cite{kac90,macd72}.
Although we have no need for the theory of affine root systems here,
it may be instructive to explain the connection to the root system 
$\mathrm{C}_n^{(1)}$ (the case of $\mathrm{A}_n^{(1)}$ is similar).
The  Weyl denominator of an affine root system $R$ is the formal
product $\prod_{\alpha\in R_+}(1-\eup^{-\alpha})^{m(\alpha)}$, where
$R_+$ denotes the set of positive roots and $m$ is a multiplicity function.
For $\mathrm{C}_n^{(1)}$, the positive roots are
\begin{align*}
m\,\delta,&  \qquad m\geq 1, \\
m\,\delta+2\varepsilon_i,& \qquad m\geq 0,\ 1\leq i\leq n, \\
m\,\delta-2\varepsilon_i,&\qquad m\geq 1,\ 1\leq i\leq n, \\
m\,\delta+\varepsilon_i\pm \varepsilon_j,& \qquad m\geq 0,\ 
1\leq i<j\leq n, \\
m\,\delta-\varepsilon_i\pm \varepsilon_j,&\qquad  m\geq 1,\
1\leq i<j\leq n, 
\end{align*}
where $\varepsilon_1,\dots,\varepsilon_n$ are the coordinate functions 
on $\R^n$ and $\delta$ is the constant function $1$.
The roots $m\,\delta$ have multiplicity $n$, while all other roots have 
multiplicity $1$.
Thus, the Weyl denominator for $\mathrm{C}_n^{(1)}$ is
\begin{multline*}
\prod_{m=0}^{\infty}\bigg(\big(1-\eup^{-(m+1)\delta}\big)^n
\prod_{i=1}^n \big(1-\eup^{-m\,\delta-2\varepsilon_i}\big)
\big(1-\eup^{-(m+1)\delta+2\varepsilon_i}\big)\\
\times\prod_{1\leq i<j\leq n}
\big(1-\eup^{-m\,\delta-\varepsilon_i-\varepsilon_j}\big)
\big(1-\eup^{-(m+1)\delta+\varepsilon_i+\varepsilon_j}\big)
\big(1-\eup^{-m\,\delta-\varepsilon_i+\varepsilon_j}\big)
\big(1-\eup^{-(m+1)\delta+\varepsilon_i-\varepsilon_j}\big)\bigg).
\end{multline*} 
It is easy to check that this equals 
$(p;p)_{\infty}^n\,x_1^0x_2^{-1}\cdots x_n^{1-n}\,
\Delta^{\mathrm{C}}(x_1,\dots,x_n)$, where 
$p=\eup^{-\delta}$ and $x_i=\eup^{-\varepsilon_i}$.

We will consider elliptic hypergeometric series containing the factor
$\Delta^{\mathrm{A}}(xq^k)$ or $\Delta^{\mathrm{C}}(xq^k)$, where
$xq^k=\big(x_1q^{k_1},x_2q^{k_2},\dots,x_rq^{k_r}\big)$, 
the $k_i{\in\mathbb{Z}}$ being summation indices, and $r=n+1$ in the 
case of $\mathrm{A}_n$ and $r=n$ in the case of $\mathrm{C}_n$. 
We refer to these as $\mathrm{A}_n$ and $\mathrm{C}_n$ series, 
respectively.  
In the case of $\mathrm{A}_n$,  the summation variables typically 
satisfy a restriction of the form $k_1+\dots+k_{n+1}=N$.
Eliminating $k_{n+1}$ gives series containing the $\mathrm{A}_{n-1}^{(1)}$ 
Weyl denominator times $\prod_{i=1}^n\theta(ax_iq^{k_i+\abs{k}})$, where
$a=q^{-N}/x_{n+1}$; these will also be viewed as $\mathrm{A}_n$ series.

Similarly, $\mathrm{A}_n$ integrals contain the factor
\begin{subequations}\label{if}
\begin{equation}\label{af}
\frac{1}{\prod_{1\leq i<j\leq n+1}\Gamma(z_i/z_j,z_j/z_i)}, 
\end{equation}
where $z_1\cdots z_{n+1}=1$, while $\mathrm{C}_n$ integrals contain
\begin{equation}\label{cf}
\frac{1}{\prod_{i=1}^n\Gamma\big(z_i^{\pm 2}\big)\prod_{1\leq i<j\leq n}
\Gamma(z_i^{\pm}z_j^{\pm})}.
\end{equation}
\end{subequations}
If we denote the expression \eqref{cf} by $g(z)$ then it is easy to 
verify that, for $k\in\mathbb{Z}^n$,
\[
\frac{g(zq^k)}{g(z)}
=\Bigg( \prod_{i=1}^nq^{-nk_i-(n+1)k_i^2}z_i^{-2(n+1)k_i}\Bigg)
\frac{\Delta^\mathrm{C}(zq^k)}{\Delta^\mathrm{C}(z)}.
\]
A similar relation holds for the $\mathrm{A}$-type factors.
This shows that the series can be considered as discrete analogues of 
the integrals.
In fact in many instances the series can be obtained from the integrals
via residue calculus.

It is customary to attach a ``type'' to hypergeometric integrals associated
with root systems, although different authors have used slightly different
definitions of type.
As the terminology will be used here, in type I integrals the only factors
containing more than one integration variable are \eqref{if}, while type
II integrals contain twice the number of such factors. 
For example, $\mathrm{C}_n^{(\text{II})}$ integrals contain the factor 
$\prod_{i<j}\Gamma(tz_i^{\pm}z_j^{\pm})/\Gamma(z_i^{\pm}z_j^{\pm})$. 
It may be noted that, under appropriate assumptions on the parameters,
\begin{align*}
\lim_{q\rightarrow 1}\lim_{p\rightarrow 0}
\prod_{i<j}\frac{\Gamma(q^{t\pm z_i\pm z_j})}{\Gamma(q^{\pm z_i\pm z_j})}
&=\lim_{q\rightarrow 1}\prod_{i<j}
\frac{(q^{\pm z_i\pm z_j};q)_{\infty}}{(q^{t\pm z_i\pm z_j};q)_{\infty}}\\
&=\prod_{i<j}(1-z_i^{\pm}z_j^{\pm})^t=
\prod_{i<j}\Big((z_i+z_i^{-1})-(z_j+z_j^{-1})\Big)^{2t}.
\end{align*}
For this reason $\mathrm{C}_n$ beta integrals of type II are sometimes
referred to as elliptic Selberg integrals.
There are also integrals containing an intermediate number of factors. 
We will refer to these as integrals of mixed type.

\section{Integrals}\label{Sec_Int}

Throughout this section we assume that $\abs{q}<1$. 
Whenever possible, we have restricted the parameters in such a way that the 
integrals may be taken over the $n$-dimensional complex torus $\mathbb{T}^n$.
However, all results can be extended to more general parameter domains
by appropriately deforming $\mathbb{T}^n$.

When $n=1$ all the stated $\mathrm{A}_n$ and $\mathrm{C}_n$ beta integral
evaluations reduce to Spiridonov's elliptic beta integral.

\subsection{$\mathrm{A}_n$ beta integrals}

We will present four $\mathrm{A}_n$ beta integrals.
In each of these the integrand contains a variable $z_{n+1}$ which is
determined from the integration variables $z_1,\dots,z_n$ by the relation 
$z_1\cdots z_{n+1}=1$. 
To shorten the expressions we define the constant $\kappa_n^{\mathrm{A}}$ 
by
\[
\kappa_n^{\mathrm{A}}=\frac{(p;p)_{\infty}^n(q;q)_{\infty}^n}
{(n+1)!(2\pi\hspace{1pt}\textup{i})^n}.
\]

For $1\leq i\leq n+2$, let $\abs{s_i}<1$ and $\abs{t_i}<1$, such that
$ST=pq$, where $S=s_1\cdots s_{n+2}$ and $T=t_1\cdots t_{n+2}$. 
Then we have the type I integral
\begin{equation}\label{ai1}
\kappa_n^{\mathrm{A}}\int_{\mathbb T^n}
\frac{\prod_{i=1}^{n+2}\prod_{j=1}^{n+1}\Gamma(s_iz_j,t_i/z_j)}
{\prod_{1\leq i<j\leq n+1}{\Gamma(z_i/z_j,z_j/z_i)}}\,
\frac{\dup z_1}{z_1}\cdots\frac{\dup z_n}{z_n} 
=\prod_{i=1}^{n+2}\Gamma(S/s_i,T/t_i)
\prod_{i,j=1}^{n+2}\Gamma(s_it_j).
\end{equation}

Next, let $\abs{s}<1$, $\abs{t}<1$, $\abs{s_i}<1$ and $\abs{t_i}<1$ 
for $1\leq i\leq 3$, such that $s^{n-1}t^{n-1}s_1s_2s_3t_1t_2t_3=pq$.
Then we have the type II integral
\begin{multline}\label{ai2}
\kappa_n^{\mathrm{A}} \int_{\mathbb T^n}
\prod_{1\leq i<j\leq n+1} \frac{\Gamma(sz_iz_j,t/z_iz_j)}
{\Gamma(z_i/z_j,z_j/z_i)}
\prod_{i=1}^3 \prod_{j=1}^{n+1}\Gamma(s_iz_j,t_i/z_j)\,
\frac{\dup z_1}{z_1}\cdots\frac{\dup z_n}{z_n}\\
=\begin{cases}
\displaystyle
\prod_{m=1}^{N} \bigg(\Gamma\big(s^mt^m\big)\prod_{1\leq i<j\leq 3}
\Gamma\big(s^{m-1}t^ms_is_j,s^mt^{m-1}t_it_j\big)
\prod_{i,j=1}^3\Gamma\big(s^{m-1}t^{m-1}s_it_j\big)\bigg)\\
\displaystyle
\qquad\times\,\Gamma\big(s^{N-1}s_1s_2s_3,t^{N-1}t_1t_2t_3\big)
\prod_{i=1}^3\Gamma\big(s^Ns_i,t^Nt_i\big), &  n=2N,\\
\displaystyle\prod_{m=1}^N
\bigg(\Gamma(s^mt^m)\prod_{1\leq i<j\leq 3}
\Gamma\big(s^{m-1}t^ms_is_j,s^mt^{m-1}t_it_j\big)\bigg)
\prod_{m=1}^{N+1}\prod_{i,j=1}^3\Gamma\big(s^{m-1}t^{m-1}s_it_j\big)\\[5mm]
\displaystyle
\qquad\times\,\Gamma\big(s^{N+1},t^{N+1}\big)
\prod_{1\leq i<j\leq 3}
\Gamma\big(s^Ns_is_j,t^Nt_it_j\big), &  n=2N+1.
\end{cases}
\end{multline}

Let $\abs{t}<1$, $\abs{t_i}<1$ for $1\leq i\leq n+3$ and 
$\abs{t}<\abs{s_i}<\abs{t}^{-1}$ for $1\leq i\leq n$, 
where $t^2t_1\cdots t_{n+3}=pq$.
Then,
\begin{multline}\label{ai3}
\kappa_n^{\mathrm{A}}\int_{\mathbb T^n}\prod_{1\leq i<j\leq n+1}
\frac{1}{\Gamma(z_i/z_j,z_j/z_i,t^2z_iz_j)}
\prod_{j=1}^{n+1} \bigg(\prod_{i=1}^n \Gamma(ts_i^{\pm} z_j)
\prod_{i=1}^{n+3}\Gamma(t_i/z_j) \bigg) \, 
\frac{\dup z_1}{z_1}\cdots\frac{\dup z_n}{z_n}  \\
=\prod_{i=1}^n\prod_{j=1}^{n+3} \Gamma(ts_i^{\pm}t_j)
\prod_{1\leq i<j\leq n+3} \frac{1}{\Gamma(t^2t_it_j)},
\end{multline}
which is an integral of mixed type.

Finally, let $\abs{t}<1$, $\abs{s_i}<1$ for $1\leq i\leq 4$ and
$\abs{t_i}<1$ for $1\leq i\leq n+1$ such that $t^{n-1}s_1\cdots s_4T=pq$, 
where $T=t_1\cdots t_{n+1}$.
Then we have a second mixed-type integral:
\begin{multline}\label{ai4}
\kappa_n^{\mathrm{A}} \int_{\mathbb T^n}
\prod_{1\leq i<j\leq n+1} \frac{\Gamma(tz_iz_j)}
{\Gamma(z_i/z_j,z_j/z_i)}\prod_{j=1}^{n+1}
\bigg(\prod_{i=1}^4\Gamma(s_iz_j)\prod_{i=1}^{n+1}
\Gamma(t_i/z_j)\bigg)\,
\frac{\dup z_1}{z_1}\cdots\frac{\dup z_n}{z_n}\\
=\begin{cases} \displaystyle
\Gamma(T)\prod_{i=1}^4\frac{\Gamma(t^Ns_i)}
{\Gamma(t^NTs_i)}\prod_{1\leq i<j\leq n+1}\Gamma(tt_it_j)
\prod_{i=1}^4 \prod_{j=1}^{n+1}\Gamma(s_it_j), & n=2N,\\
\displaystyle
\frac{\Gamma(t^{N+1},T)}{\Gamma(t^{N+1}T)}
\prod_{1\leq i<j\leq 4}\Gamma(t^Ns_is_j)
\prod_{1\leq i<j\leq n+1}\Gamma(tt_it_j)
\prod_{i=1}^4 \prod_{j=1}^{n+1} \Gamma(s_it_j),& n=2N+1.
\end{cases}
\end{multline}

\subsection{$\mathrm{C}_n$ beta integrals}

We will give three $\mathrm{C}_n$ beta integrals. 
They all involve the constant
\[
\kappa_n^{\mathrm{C}}=\frac{(p;p)_{\infty}^n(q;q)_{\infty}^n}
{n!2^n(2\pi\hspace{1pt}\iup)^n}.
\]

Let $\abs{t_i}<1$ for $1\leq i\leq 2n+4$ such that $t_1\cdots t_{2n+4}=pq$. 
We then have the following $\mathrm{C}_n$ beta integral of type I
\begin{equation}\label{ci1}
\kappa_n^{\mathrm{C}} \int_{\mathbb T^n}
\prod_{1\leq i<j\leq n} \frac{1}{\Gamma(z_i^{\pm}z_j^{\pm})}
\prod_{j=1}^n\frac{\prod_{i=1}^{2n+4}\Gamma(t_iz_j^{\pm})}
{\Gamma(z_j^{\pm 2})}\,
\frac{\dup z_1}{z_1}\cdots\frac{\dup z_n}{z_n} 
=\prod_{1\leq i<j\leq 2n+4}\Gamma(t_it_j).
\end{equation}

Next, let $\abs{t}<1$ and $\abs{t_i}<1$ for $1\leq i\leq 6$ such that
$t^{2n-2}t_1\cdots t_6=pq$.
We then have the type II $\mathrm{C}_n$ beta integral 
\begin{equation}\label{ci2}
\kappa_n^{\mathrm{C}} \int_{\mathbb T^n}
\prod_{1\leq i<j\leq n} \frac{\Gamma(tz_i^{\pm}z_j^{\pm})}
{\Gamma(z_i^{\pm}z_j^{\pm})}
\prod_{j=1}^n\frac{\prod_{i=1}^6\Gamma(t_iz_j^{\pm})}
{\Gamma(z_j^{\pm 2})}\,
\frac{\dup z_1}{z_1}\cdots\frac{\dup z_n}{z_n} 
=\prod_{m=1}^n\bigg(\,\frac{\Gamma(t^m)}{\Gamma(t)}
\prod_{1\leq i<j\leq 6}\Gamma(t^{m-1}t_it_j)\bigg).
\end{equation}
This is the elliptic Selberg integral mentioned in the introduction.

At this point it is convenient to introduce notation for more general
$\mathrm{C}_n$ integrals of type II. 
For $m$ a non-negative integer, let $\abs{t}<1$ and $\abs{t_i}<1$ 
for $1\leq i\leq 2m+6$ such that
\begin{equation}\label{c2b}
t^{2n-2}t_1\cdots t_{2m+6}=(pq)^{m+1}.
\end{equation} 
We then define
\begin{equation}\label{c2not}
J_{\mathrm{C}_n}^{(m)}(t_1,\dots,t_{2m+6};t)=
\kappa_n^{\mathrm{C}} \int_{\mathbb T^n}
\prod_{1\leq i<j\leq n} \frac{\Gamma(tz_i^{\pm}z_j^{\pm})}
{\Gamma(z_i^{\pm}z_j^{\pm})}
\prod_{j=1}^n\frac{\prod_{i=1}^{2m+6}\Gamma(t_iz_j^{\pm})}
{\Gamma(z_j^{\pm 2})}\,
\frac{\dup z_1}{z_1}\cdots\frac{\dup z_n}{z_n}. 
\end{equation}
Note that \eqref{ci2} gives a closed-form evaluation for
the integral $J_{\mathrm{C}_n}^{(0)}$. 
As outlined in \cite[Appendix]{rains10}, $J_{\mathrm{C}_n}^{(m)}$ can be 
continued to a single-valued meromorphic function in the parameters $t_i$
and $t$ subject to the constraint \eqref{c2b}. 
For generic values of the parameters this continuation is obtained 
by replacing the integration domain with an appropriate deformation of
$\mathbb T^n$.
We can now state the second $\mathrm{C}_n$ beta integral of type II as
\begin{multline}\label{ci3}
J_{\mathrm{C}_n}^{(n-1)}(t_1,\dots,t_4,s_1,\dots,s_n,pq/ts_1,\dots,pq/ts_n;t)
\\
=\Gamma(t)^n
\prod_{l=1}^n\prod_{1\leq i<j\leq 4}\Gamma\big(t^{l-1}t_it_j\big)
\prod_{i=1}^n\prod_{j=1}^4 \frac{\Gamma(s_it_j)}{\Gamma(ts_i/t_j)},
\end{multline}
where  $t^{n-2}t_1t_2t_3t_4=1$. 
In this identity it is necessary to work with an analytic continuation of
\eqref{c2not} since the inequalities $\abs{t_i},\,\abs{t}<1$ are 
incompatible with $t^{n-2}t_1t_2t_3t_4=1$ for $n\geq 2$.

\subsection{Integral transformations}

We  now turn to integral transformations, starting with integrals of type I. 
For $m$ a non-negative integer we introduce the notation
\[
I_{\mathrm{A}_n}^{(m)}(s_1,\dots,s_{m+n+2};t_1,\dots,t_{m+n+2})=
\kappa_n^{\mathrm{A}}\int_{\mathbb T^n}\frac{
\prod_{i=1}^{m+n+2}\prod_{j=1}^{n+1}\Gamma(s_iz_j,t_i/z_j)}
{\prod_{1\leq i<j\leq n+1}{\Gamma(z_i/z_j,z_j/z_i)}}\,
\frac{\dup z_1}{z_1}\cdots\frac{\dup z_n}{z_n},
\]
where $\abs{s_i}<1$ and $\abs{t_i}<1 $ for all $i$, 
$\prod_{i=1}^{m+n+2}s_it_i=(pq)^{m+1}$ and $z_1\cdots z_{n+1}=1$. 
We also define
\[
I_{\mathrm{C}_n}^{(m)}(t_1,\dots,t_{2m+2n+4})
=\kappa_n^{\mathrm{C}} \int_{\mathbb T^n}
\prod_{1\leq i<j\leq n} \frac{1}{\Gamma(z_i^{\pm}z_j^{\pm})}
\prod_{j=1}^n\frac{\prod_{i=1}^{2m+2n+4}\Gamma(t_iz_j^{\pm})}
{\Gamma(z_j^{\pm 2})}\,
\frac{\dup z_1}{z_1}\cdots\frac{\dup z_n}{z_n},
\]
where  $\abs{t_i}<1$ for all $i$ and $t_1\cdots t_{2m+2n+4}=(pq)^{m+1}$. 
The $\mathrm{A}_n$ integral satisfies
\[
I_{\mathrm{A}_n}^{(m)}(s_1,\dots,s_{m+n+2};t_1,\dots,t_{m+n+2})=
I_{\mathrm{A}_n}^{(m)}(s_1\zeta,\dots,s_{m+n+2}\zeta;
t_1/\zeta,\dots,t_{m+n+2}/\zeta)
\]
for $\zeta$ any $(n+1)$-th root of unity, whereas the $\mathrm{C}_n$
integral is invariant under simultaneous negation of all of the $t_i$.
We further note that \eqref{ai1} and \eqref{ci1} provide closed-form 
evaluations of $I_{\mathrm{A}_n}^{(0)}$ and $I_{\mathrm{C}_n}^{(0)}$,
respectively.

For the integral $I_{\mathrm{A}_n}^{(m)}$, the following transformation
reverses the roles of $m$ and $n$:
\begin{multline}\label{ait}
I_{\mathrm{A}_n}^{(m)}(s_1,\dots,s_{m+n+2};t_1,\dots,t_{m+n+2})\\
=\prod_{i,j=1}^{m+n+2}\Gamma(s_it_j)\cdot I_{\mathrm{A}_m}^{(n)}
\Bigg(\frac{\lambda}{s_1},\dots,\frac{\lambda}
{s_{m+n+2}};\frac{pq}{\lambda t_1},\dots,
\frac{pq}{\lambda t_{m+n+2}}\Bigg),
\end{multline}
where $\lambda^{m+1}=s_1\cdots s_{m+n+2}$, 
$(pq/\lambda)^{m+1}=t_1\cdots t_{m+n+2}$.
Moreover, for $t_1\cdots t_{2m+2n+4}=(pq)^{m+1}$, there is an analogous 
transformation of type $\mathrm{C}$:
\begin{equation}\label{cit}
I_{\mathrm{C}_n}^{(m)}(t_1,\dots,t_{2m+2n+4})
=\prod_{1\leq i<j\leq 2m+2n+4}\Gamma(t_it_j)\cdot I_{\mathrm{C}_m}^{(n)}
\Bigg(\frac{\sqrt{pq}}{t_1},\dots,\frac{\sqrt{pq}}{t_{2m+2n+4}}\Bigg).
\end{equation}

It is easy to check that 
$I_{\mathrm{A}_1}^{(m)}(t_1,\dots,t_{m+3};t_{m+4},\dots,t_{2m+6})=
I_{\mathrm{C}_1}^{(m)}(t_1,\dots,t_{2m+6})$.
Thus, combining \eqref{ait} and \eqref{cit} leads to
\begin{multline}\label{acit}
I_{\mathrm{A}_n}^{(1)}(s_1,\dots,s_{n+3};t_1,\dots,t_{n+3}) \\
=\prod_{1\leq i<j\leq n+3}
\Gamma(S/s_is_j,T/t_it_j)\cdot I_{\mathrm{C}_n}^{(1)}
(s_1/v,\dots,s_{n+3}/v,t_1v,\dots,t_{n+3}v),
\end{multline}
where $S=s_1\cdots s_{n+3}$, $T=t_1\cdots t_{n+3}$ and $\nu^2=S/pq=pq/T$.
Since $I_{\mathrm{C}_n}^{(1)}$ is symmetric,  \eqref{acit} implies
non-trivial symmetries of $I_{\mathrm{A}_n}^{(1)}$, such as
\begin{multline}\label{ait2}
I_{\mathrm{A}_n}^{(1)}(s_1,\dots,s_{n+3};t_1,\dots,t_{n+3})=
\prod_{i=1}^{n+2}\Gamma(s_it_{n+3},t_is_{n+3},S/s_is_{n+3},T/t_it_{n+3})\\
\times I_{\mathrm{A}_n}^{(1)}(s_1/v,\dots,s_{n+2}/v,s_{n+3}v^n;
t_1v,\dots,t_{n+2}v,t_{n+3}/v^n),
\end{multline}
where, with the same definitions of $S$ and $T$ as above, $ST=(pq)^2$
and $\nu^{n+1}=St_{n+3}/pqs_{n+3}$.

\medskip

\begin{multline}\label{cit2}
J_{\mathrm{C}_n}^{(1)}(t_1,\dots,t_8;t)=
\prod_{m=1}^n \bigg(\prod_{1\leq i<j\leq 4} \Gamma(t^{m-1}t_it_j)
\prod_{5\leq i<j\leq 8} \Gamma(t^{m-1}t_it_j)\bigg) \\
\times
J_{\mathrm{C}_n}^{(1)}(t_1v,\dots,t_4v,t_5/v,\dots,t_8/v;t),
\end{multline}
where 
$t^{2n-2}t_1\cdots t_8=(pq)^2$ and 
$v^2=pqt^{1-n}/t_1t_2t_3t_4=t^{n-1}t_5t_6t_7t_8/pq$.
Iterating this transformation yields a symmetry of
$J_{\mathrm{C}_n}^{(1)}$ under the Weyl group of type $\mathrm{E}_7$ 
\cite{rains10}.

We conclude with a transformation between $\mathrm{C}_n$ and $\mathrm{C}_m$
integrals of type II:
\begin{align}
J_{\mathrm{C}_n}^{(m+n-1)}&(t_1,\dots,t_4,s_1,\dots,
s_{m+n},pq/ts_1,\dots,pq/ts_{m+n};t) \notag \\
&=\Gamma(t)^{n-m}
\prod_{1\leq i<j\leq 4} \frac{\prod_{l=1}^n \Gamma\big(t^{l-1}t_it_j\big)}
{\prod_{l=1}^m \Gamma\big(t^{l+n-m-1}t_it_j\big)}
\prod_{i=1}^{m+n} \prod_{j=1}^4
\frac{\Gamma(s_it_j)}{\Gamma(ts_i/t_j)} \notag \\
&\qquad\times J_{\mathrm{C}_m}^{(m+n-1)}(t/t_1,\dots,t/t_4,s_1,\dots,
s_{m+n},pq/ts_1,\dots,pq/ts_{m+n};t),
\label{cit3}
\end{align}
where $t_1t_2t_3t_4=t^{m-n+2}$.

\subsection{Notes}\label{Sec_intnotes}

For $p=0$ the integrals \eqref{ai1}, \eqref{ai2}, \eqref{ci1} and 
\eqref{ci2} are due to Gustafson \cite{gust92, gust94}, the integral
\eqref{ai4} to Gustafson and Rakha \cite{gustrakha00} and the
transformation \eqref{ait2} to Denis and Gustafson \cite{denisgust92}. 
None of the $p=0$ instances of \eqref{ai3}, \eqref{ci3}--\eqref{acit}, 
\eqref{cit2} and \eqref{cit3} were known prior to the elliptic case.

For general $p$, van Diejen and Spiridonov conjectured the type I 
$\mathrm{C}_n$ beta integral \eqref{ci1} and showed that it implies the
elliptic Selberg integral \eqref{ci2} \cite{vandspir00,vandspir01}.
A rigorous derivation of the classical Selberg integral as a special
limit of \eqref{ci2} is due to Rains \cite{rains09}.
Spiridonov \cite{spir04} conjectured the type I $\mathrm{A}_n$ beta
integral \eqref{ai1} and showed that, combined with \eqref{ci1}, it
implies the type II $\mathrm{A}_n$ beta integral \eqref{ai2}, as well
as the integral \eqref{ai4} of mixed type.
He also showed that \eqref{ai1} implies \eqref{ait2}. 
The first proofs of the fundamental type I integrals \eqref{ai1} and 
\eqref{ci1} were obtained by Rains \cite{rains10}.
For subsequent proofs of \eqref{ai1}, \eqref{ci1} and \eqref{ci2},
see \cite{spir07}, \cite{rainsspir09,spir07} and \cite{itonoumi17b}, 
respectively.
In \cite{rains10} Rains also proved the integral transformations 
\eqref{ait}, \eqref{cit} and \eqref{cit2}, and gave further 
transformations analogous to \eqref{cit2}.
The integral \eqref{ai3} of mixed type is due to Spiridonov and
Warnaar \cite{spirwar06}.
The transformation \eqref{cit3}, which includes \eqref{ci3} as its
$m=0$ case, was conjectured by Rains \cite{rains12} and also appears 
in \cite{spirvart11}.
It was first proved by van der Bult in \cite{vdb09} and subsequently
proved and generalized to an identity for the ``interpolation kernel''
(an analytic continuation of the elliptic interpolation functions
$R^{\ast}_{\lambda}$ of Section~\ref{Sec_EMK}) in \cite{rains18}.

Several of the integral identities surveyed here have analogues for
$\abs{q}=1$.
In the case of  \eqref{ai1}, \eqref{ai3} and \eqref{ci1} these were found
in \cite{vandspir05}, and the unit-circle analogue of \eqref{ci2} is given
in \cite{vandspir05}.

In \cite{spir04} Spiridonov gives one more $\mathrm{C}_n$ beta integral, 
which lacks the $p\leftrightarrow q$ symmetry present in all the integrals 
considered here, and is more elementary in that it follows as a determinant
of one-variable beta integrals.

In \cite{rains12} Rains conjectured several quadratic integral
transformations involving the interpolation functions
$R^{*}_{\lambda}$. 
These conjectures were proved in \cite{vdb11,rains18}. 
In special cases, they simplify to transformations for the function 
$J_{\mathrm{C}_n}^{(2)}$.

Motivated by quantum field theories on lens spaces, Spiridonov 
\cite{spir18} evaluated certain finite sums of $\mathrm{C}_n$ integrals,
both for type I and type II.
In closely related work, Kels and Yamazaki \cite{kels18} obtained 
transformation formulas for finite sums of $\mathrm{A}_n$ and $\mathrm{C}_n$
integrals of type I.

As mentioned in the introduction, the recent identification of elliptic 
hypergeometric integrals as indices in supersymmetric quantum field theory
by Dolan and Osborn \cite{dolosb09} has led to a large number of conjectured
integral evaluations and transformations
\cite{gadde10,gadde10b,spirvart11,spirvart12,spirvart14}. 
It is too early to give a survey of the emerging picture, but it is
clear that the identities stated in this section are a small sample 
from a much larger collection of identities.

\section{Series}\label{Sec_Series}

In this section we give the most important summation and transformation
formulas for elliptic hypergeometric series associated to $\mathrm{A}_n$
and $\mathrm{C}_n$. 
In the $n=1$ case all summations except for \eqref{bs} simplify to the 
elliptic Jackson summation of Frenkel and Turaev.
Similarly, most transformations may be viewed as generalizations of the
elliptic Bailey transformation.

\subsection{$\mathrm{A}_n$ summations}

The following $\mathrm{A}_n$ elliptic Jackson summation is a discrete 
analogue of the multiple beta integral \eqref{ai1}: 
\begin{subequations}\label{as1}
\begin{equation}\label{as1a}
\sum_{\substack{k_1,\dots,k_{n+1}\geq 0\\ k_1+\cdots+k_{n+1}=N}}
\frac{\Delta^{\mathrm{A}}(xq^k)}{\Delta^{\mathrm{A}}(x)}
\prod_{i=1}^{n+1}\frac{\prod_{j=1}^{n+2}(x_ia_j)_{k_i}}
{(bx_i)_{k_i}\prod_{j=1}^{n+1}(qx_i/x_j)_{k_i}}
=\frac{(b/a_1,\dots,b/a_{n+2})_N}
{(q,bx_1,\dots,bx_{n+1})_N},
\end{equation}
where $b=a_1\dotsm a_{n+2}x_1\cdots x_{n+1}$.
Using the constraint on the summation indices to eliminate $k_{n+1}$, 
this identity can be written less symmetrically as
\begin{multline}\label{as1b}
\sum_{\substack{k_1,\dots, k_n\geq 0\\ \abs{k}\leq N}}
\frac{\Delta^{\mathrm{A}}(xq^k)}{\Delta^{\mathrm{A}}(x)}
\prod_{i=1}^n\Bigg(\frac{\theta(ax_iq^{k_i+\abs{k}})}{\theta(ax_i)}
\,\frac{(ax_i)_{\abs{k}}\prod_{j=1}^{n+2}(x_ib_j)_{k_i}}
{(aq^{N+1}x_i,aqx_i/c)_{k_i}\prod_{j=1}^n(qx_i/x_j)_{k_i}}\Bigg) 
\frac{(q^{-N},c)_{\abs{k}}}
{\prod_{i=1}^{n+2}(aq/b_i)_{\abs{k}}}\,q^{\abs{k}} \\
=c^N\prod_{i=1}^n\frac{(aqx_i)_N}{(aqx_i/c)_N}\prod_{i=1}^{n+2}
\frac{(aq/cb_i)_N}{(aq/b_i)_N},
\end{multline}
where $b_1\dotsm b_{n+2}cx_1\dotsm x_n=a^2q^{N+1}$.
By analytic continuation one can then deduce the companion identity
\begin{align}
\sum_{k_1,\dots,k_n=0}^{N_1,\dots,N_n}&\Bigg(
\frac{\Delta^{\mathrm{A}}(xq^k)}{\Delta^{\mathrm{A}}(x)}
\prod_{i=1}^n\bigg(\frac{\theta(ax_iq^{k_i+\abs{k}})}{\theta(ax_i)}\,
\frac{(ax_i)_{\abs{k}}(dx_i,ex_i)_{k_i}}
{(aq^{N_i+1}x_i)_{\abs{k}}(aqx_i/b,aqx_i/c)_{k_i}}\bigg) \notag \\
& \times\frac{(b,c)_{\abs{k}}\, q^{\abs{k}}}{(aq/d,aq/e)_{\abs{k}}}\,
\prod_{i,j=1}^n\frac{(q^{-N_j}x_i/x_j)_{k_i}}{(qx_i/x_j)_{k_i}} \Bigg)
=\frac{(aq/cd,aq/bd)_{\abs{N}}}{(aq/d,aq/bcd)_{\abs{N}}}
\prod_{i=1}^n\frac{(aqx_i,aqx_i/bc)_{N_i}}{(aqx_i/b,aqx_i/c)_{N_i}},
\label{as1c}
\end{align}
where $bcde=a^2q^{\abs{N}+1}$.
\end{subequations}

For the discrete analogue of the type II integral \eqref{ai2} we refer
the reader to the Notes at the end of this section.

Our next result corresponds to a discretization of \eqref{ai3}:
\begin{multline}\label{as2}
\sum_{\substack{k_1,\dots,k_{n+1}\geq 0 \\ k_1+\dots+k_{n+1}=N}}
\frac{\Delta^{\mathrm{A}}(xq^k)}{\Delta^{\mathrm{A}}(x)}
\prod_{1\leq i<j\leq n+1}\frac{1}{(x_ix_j)_{k_i+k_j}}
\prod_{i=1}^{n+1}\frac{q^{\binom{k_i}2}x_i^{k_i}\prod_{j=1}^n
(x_ia_j^{\pm})_{k_i}}{(bx_i,q^{1-N}x_i/b)_{k_i}
\prod_{j=1}^{n+1}(qx_i/x_j)_{k_i}} \\
=\big({-}bq^{N-1}\big)^N\frac{\prod_{i=1}^n(ba_i^{\pm})_N}
{(q)_N\prod_{i=1}^{n+1}(bx_i^{\pm})_N}.
\end{multline}
Mimicking the steps that led from \eqref{as1a} to \eqref{as1c}, the
identity \eqref{as2} can be rewritten as a sum over an $n$-dimensional
rectangle, see \cite{ros04}.
Some authors have associated \eqref{as2} and related results with the
root system $\mathrm{D}_n$ rather than $\mathrm{A}_n$.

Finally, the following summation is a discrete analogue of \eqref{ai4}:
\begin{multline}\label{as4}
\sum_{\substack{k_1,\dots,k_{n+1}\geq 0,\\k_1+\dots+k_{n+1}=N}} 
\frac{\Delta^A(xq^k)}{\Delta^A(x)}
\prod_{1\leq i<j\leq n+1}q^{k_ik_j}(x_ix_j)_{k_i+k_j}\prod_{i=1}^{n+1}
\frac{\prod_{j=1}^4(x_ib_j)_{k_i}}
{x_i^{k_i}\prod_{j=1}^{n+1}(qx_i/x_j)_{k_i}}\\
=\begin{cases}
\displaystyle\frac{(Xb_1,Xb_2,Xb_3,Xb_4)_N}{X^N(q)_N}, 
& \text{$n$ odd},\\[4mm]
\displaystyle\frac{(X,Xb_1b_2,Xb_1b_3,Xb_1b_4)_N}{(Xb_1)^N(q)_N}, 
& \text{$n$ even},
\end{cases}
\end{multline}
where $X=x_1\dotsm x_{n+1}$ and $q^{N-1}b_1\dotsm b_4X^2=1$.

\subsection{$\mathrm{C}_n$ summations}

The following $\mathrm{C}_n$ elliptic Jackson summation is a discrete 
analogue of \eqref{ci1}:
\begin{multline}\label{dgs}
\sum_{k_1,\dots,k_n=0}^{N_1,\dots,N_n}
\frac{\Delta^{\mathrm{C}}(xq^k)}{\Delta^{\mathrm{C}}(x)}
\prod_{i=1}^n\frac{(bx_i,cx_i,dx_i,ex_i)_{k_i}\,q^{k_i}}
{(qx_i/b,qx_i/c,qx_i/d,qx_i/e)_{k_i}} 
\prod_{i,j=1}^n\frac{(q^{-N_j}x_i/x_j,x_ix_j)_{k_i}}
{(qx_i/x_j,q^{{N_j+1}}x_ix_j)_{k_i}}\\
=\frac{\prod_{i,j=1}^n(qx_ix_j)_{N_i}}
{\prod_{1\leq i<j\leq n}(qx_ix_j)_{N_i+N_j}}
\,\frac{(q/bc,q/bd,q/cd)_{\abs{N}}}
{\prod_{i=1}^n(qx_i/b,qx_i/c,qx_i/d,q^{-N_i}e/x_i)_{N_i}},
\end{multline}
where $bcde=q^{\abs{N}+1}$.

The discrete analogues of the type II integrals \eqref{ci2} and \eqref{ci3} 
are most conveniently expressed in terms of the series
\begin{multline}\label{ds}
\V{r+1}{r}{n}(a;b_1,\dots,b_{r-4})
=\sum_{\lambda}\Bigg(\prod_{i=1}^n 
\frac{\theta(at^{2-2i}q^{2\lambda_i})}{\theta(at^{2-2i})} \,
\frac{(at^{1-n},b_1,\dots,b_{r-4})_{\lambda}}
{(qt^{n-1},aq/b_1,\dots,aq/b_{r-4})_{\lambda}}\\
\times\prod_{1\leq i<j\leq n} \bigg(
\frac{\theta(t^{j-i}q^{\lambda_i-\lambda_j},at^{2-i-j}q^{\lambda_i+\lambda_j})}
{\theta(t^{j-i},at^{2-i-j})}\,
\frac{(t^{j-i+1})_{\lambda_i-\lambda_j}(at^{3-i-j})_{\lambda_i+\lambda_j}}
{(qt^{j-i-1})_{\lambda_i-\lambda_j}(aqt^{1-i-j})_{\lambda_i+\lambda_j}}\bigg)\,
q^{\abs{\lambda}}t^{2n(\lambda)}\Bigg),
\end{multline}
where the summation is over partitions $\lambda=(\lambda_1,\dots,\lambda_n)$ 
of length at most $n$. Note that this implicitly depends on $t$ as well as 
$p$ and $q$.
When $b_{r-4}=q^{-N}$ with $N\in\mathbb Z_{\geq 0}$, this becomes a terminating
series, with sum ranging over partitions $\lambda\subset(N^n)$.
The series \eqref{ds} is associated with $\mathrm{C}_n$ since
\[
\prod_{i=1}^n \frac{\theta(at^{2-2i}q^{2\lambda_i})}
{\theta(at^{2-2i})} \prod_{1\leq i<j\leq n}
\frac{\theta(t^{j-i} q^{\lambda_i-\lambda_j},
at^{2-i-j}q^{\lambda_i+\lambda_j})}{\theta(t^{j-i},at^{2-i-j})}
=q^{-n(\lambda)}
\frac{\Delta^{\mathrm{C}}(xq^{\lambda})}{\Delta^{\mathrm{C}}(x)},
\]
with $x_i=\sqrt at^{1-i}$.

Using the above notation, the discrete analogue of \eqref{ci2} is
\begin{equation}\label{w}
\V{10}{9}{n}\big(a;b,c,d,e,q^{-N}\big)=
\frac{(aq,aq/bc,aq/bd,aq/cd)_{(N^n)}}
{(aq/b,aq/c,aq/d,aq/bcd)_{(N^n)}},
\end{equation}
where $bcdet^{n-1}=a^2q^{N+1}$. 
This is the $\mathrm{C}_n$ summation
mentioned in the introduction.

Next, we give a discrete analogue of \eqref{ci3}:
\begin{multline}\label{bs}
\V{2r+8}{2r+7}{n}\bigg(a;t^{1-n}b,\frac{a}{b},c_1q^{k_1},\dots,c_rq^{k_r},
\frac{aq}{c_1},\dots,\frac{aq}{c_r },q^{-N}\bigg)\\
=\frac{(aq,qt^{n-1})_{(N^n)}}{(bq,aqt^{n-1}/b)_{(N^n)}}
\prod_{i=1}^r\frac{(c_ib/a,c_it^{n-1}/b)_{(k_i^n)}}
{(c_i,c_it^{n-1}/a)_{(k_i^n)}},
\end{multline}
where the $k_i$ are non-negative integers such that $k_1+\dots+k_r=N$.
As the summand contains the factors 
$(c_iq^{k_i})_{\lambda}/(c_i)_{\lambda}$, 
this is a so-called Karlsson--Minton-type summation.

Finally, we have the $\mathrm{C}_n$ summation 
\begin{multline}\label{ss}
\sum_{k_1,\dots,k_n=0}^N\frac{\Delta^{\mathrm{C}}(xq^k)}
{\Delta^{\mathrm{C}}(x)} \prod_{i=1}^n
\frac{(x_i^2,bx_i,cx_i,dx_i,ex_i,q^{-N})_{k_i}}
{(q,qx_i/b,qx_i/c,qx_i/d,qx_i/e,q^{N+1}x_i^2)_{k_i}} \, q^{k_i} \\
=\prod_{1\leq i<j\leq n} \frac{\theta(q^{N}x_ix_j)}
{\theta(x_ix_j)}
\prod_{i=1}^n\frac{(qx_i^2,q^{2-i}/bc,q^{2-i}/bd,q^{2-i}/cd)_N}
{(qx_i/b,qx_i/c,qx_i/d,q^{-N}e/x_i)_N},
\end{multline}
where $bcde=q^{N-n+2}$. 
There is a  corresponding integral evaluation \cite{spir04}, which was
mentioned in \S \ref{Sec_intnotes}.

\subsection{Series transformations}

Several of the transformations stated below have companion identities 
(similar to the different versions of \eqref{as1}) which will not be 
stated explicitly.

The following $\mathrm{A}_n$ Bailey transformation is a discrete analogue
of \eqref{ait2}:
\begin{align}
&\quad\sum_{k_1,\dots,k_n=0}^{N_1,\dots,N_n}
\Bigg(\frac{\Delta^{\mathrm{A}}(xq^k)}{\Delta^{\mathrm{A}}(x)}
\prod_{i=1}^n \bigg(\frac{\theta(ax_iq^{k_i+\abs{k}})}{\theta(ax_i)}\,
\frac{(ax_i)_{\abs{k}}(ex_i,fx_i,gx_i)_{k_i}}
{(aq^{{N_i+1}}x_i)_{\abs{k}}(aqx_i/b,aqx_i/c,aqx_i/d)_{k_i}}\bigg) \notag \\
&\qquad\qquad\quad\times
\frac{(b,c,d)_{\abs{k}}\, q^{\abs{k}}}
{(aq/e,aq/f,aq/g)_{\abs{k}}} \,
\prod_{i,j=1}^n \frac{(q^{-N_j}x_i/x_j)_{k_i}}{(qx_i/x_j)_{k_i}} 
\Bigg) \notag \\
&=\Big(\frac{a}{\lambda}\Big)^{\abs{N}}
\frac{(\lambda q/f,\lambda q/g)_{\abs{N}}}{(aq/f,aq/g)_{\abs{N}}}
\prod_{i=1}^n\frac{(aqx_i,\lambda qx_i/d)_{N_i}}
{(\lambda qx_i,aqx_i/d)_{N_i}}\notag \\
&\quad\times\sum_{k_1,\dots,k_n=0}^{N_1,\dots,N_n}
\Bigg(\frac{\Delta^{\mathrm{A}}(xq^k)}{\Delta^{\mathrm{A}}(x)}
\prod_{i=1}^n\bigg(\frac{\theta(\lambda x_i q^{k_i+\abs{k}})}
{\theta(\lambda x_i)}\,
\frac{(\lambda x_i)_{\abs{k}}(\lambda ex_i/a,fx_i,gx_i)_{k_i}}
{(\lambda q^{{N_i+1}}x_i)_{\abs{k}}(aqx_i/b,aqx_i/c,\lambda qx_i/d)_{k_i}}
\bigg) \notag \\
&\qquad\qquad\quad\quad\times
\frac{(\lambda b/a,\lambda c/a,d)_{\abs{k}}\, q^{\abs{k}}}
{(aq/e,\lambda q/f,\lambda q/g)_{\abs{k}}}\,
\prod_{i,j=1}^n \frac{(q^{-N_j}x_i/x_j)_{k_i}}{(qx_i/x_j)_{k_i}}
\Bigg),
\label{ast1}
\end{align}
where $bcdefg=a^3q^{\abs{N}+2}$ and $\lambda=a^2q/bce$.
For $be=aq$ the sum on the right trivializes and
the transformation simplifies to \eqref{as1c}.

The next transformation, which relates an $\mathrm{A}_n$ and  a
$\mathrm{C}_n$ series, is a discrete analogue of \eqref{acit}:
\begin{align}
&\quad \sum_{k_1,\dots,k_n=0}^{N_1,\dots,N_n}
\frac{\Delta^{\mathrm{C}}(xq^k)}{\Delta^{\mathrm{C}}(x)}
\prod_{i=1}^n\frac{(bx_i,cx_i,dx_i,ex_i,fx_i,gx_i)_{k_i}\, q^{k_i}}
{(qx_i/b,qx_i/c,qx_i/d,qx_i/e,qx_i/f,qx_i/g)_{k_i}} 
\prod_{i,j=1}^n\frac{(q^{-N_j}x_i/x_j,x_ix_j)_{k_i}}
{(qx_i/x_j,q^{{N_j+1}}x_ix_j)_{k_i}}
 \notag \\
&=\frac{\prod_{i,j=1}^n(qx_ix_j)_{N_i}}
{\prod_{1\leq i<j\leq n}(qx_ix_j)_{N_i+N_j}}\,
\frac{(\lambda q/e,\lambda q/f,q/ef)_{\abs{N}}}
{\prod_{i=1}^n(\lambda qx_i,qx_i/e,qx_i/f,q^{-N_i}g/x_i)_{N_i}}\notag \\
&\quad\times\sum_{k_1,\dots,k_n=0}^{N_1,\dots,N_n}
\Bigg(\frac{\Delta^{\mathrm{A}}(xq^k)}{\Delta^{\mathrm{A}}(x)}
\prod_{i=1}^n\bigg(\frac{\theta(\lambda x_i q^{k_i+\abs{k}})}
{\theta(\lambda x_i)}\,
\frac{(\lambda x_i)_{\abs{k}}(ex_i,fx_i,gx_i)_{k_i}}
{(\lambda q^{{N_i+1}}x_i)_{\abs{k}}(qx_i/b,qx_i/c,qx_i/d)_{k_i}}\bigg) 
\notag \\
&\qquad\qquad\quad\times
\frac{(\lambda b,\lambda c,\lambda d)_{\abs{k}}\, q^{\abs{k}}}
{(\lambda q/e,\lambda q/f,\lambda q/g)_{\abs{k}}}
\prod_{i,j=1}^n\frac{(q^{-N_j}x_i/x_j)_{k_i}}{(qx_i/x_j)_{k_i}}
\Bigg),
\label{acst}
\end{align}
where $bcdefg=q^{\abs{N}+2}$ and $\lambda=q/bcd$.
For $bc=q$ this reduces to \eqref{dgs}.

The discrete analogue of \eqref{ait} provides a duality between 
$\mathrm{A}_n$ and $\mathrm{A}_m$ elliptic hypergeometric series:
\begin{multline}\label{te}
\sum_{\substack{k_1,\dots,k_{n+1}\geq 0\\ k_1+\dots+k_{n+1}=N}}
\frac{\Delta^{\mathrm{A}}(xq^k)}{\Delta^{\mathrm{A}}(x)}
\prod_{i=1}^{n+1}\frac{\prod_{j=1}^{m+n+2}(x_ia_j)_{k_i}}
{\prod_{j=1}^{m+1}(x_iy_j)_{k_i}\prod_{j=1}^{n+1}(qx_i/x_j)_{k_i}}\\
=\sum_{\substack{k_1,\dots,k_{m+1}\geq 0\\ k_1+\dots+k_{m+1}=N}}
\frac{\Delta^{\mathrm{A}}(yq^k)}{\Delta^{\mathrm{A}}(y)}
\prod_{i=1}^{m+1}\frac{\prod_{j=1}^{m+n+2}(y_i/a_j)_{k_i}}
{\prod_{j=1}^{n+1}(y_ix_j)_{k_i}\prod_{j=1}^{m+1}(qy_i/y_j)_{k_i}},
\end{multline}
where $w_1\cdots w_{m+1}=x_1\cdots x_{n+1}a_1\cdots a_{m+n+2}$.
For $m=0$ this reduces to \eqref{as1a}.

We next give a discrete analogue of \eqref{cit}. 
When $M_i$ and $N_i$ for $i=1,\dots,n$ are non-negative integers and 
$bcde=q^{\abs{N}-\abs{M}+1}$, then
\begin{align}
&\quad\sum_{k_1,\dots,k_n=0}^{N_1,\dots,N_n}
\Bigg(\frac{\Delta^{\mathrm{C}}(xq^k)}{\Delta^{\mathrm{C}}(x)}
\prod_{i=1}^n\frac{(bx_i,cx_i,dx_i,ex_i)_{k_i}\,q^{k_i}}
{(qx_i/b,qx_i/c,qx_i/d,qx_i/e)_{k_i}} \notag \\
&\qquad\qquad\times
\prod_{i=1}^n\prod_{j=1}^m 
\frac{(q^{M_j}x_iy_j,qx_i/y_j)_{k_i}}
{(x_iy_j,q^{1-M_j}x_i/y_j)_{k_i}}
\prod_{i,j=1}^n\frac{(q^{-N_j}x_i/x_j,x_ix_j)_{k_i}}
{(qx_i/x_j,q^{{N_j+1}}x_ix_j)_{k_i}}\Bigg) \notag \\
&=q^{-\abs{N}\abs{M}} \frac{(q/bc,q/bd,q/cd)_{\abs{N}}}
{(q^{-\abs{N}}bc,q^{-\abs{N}}bd,q^{-\abs{N}}cd)_{\abs{M}}}
\prod_{i=1}^m\prod_{j=1}^n\frac{(q^{-N_j}y_i/x_j)_{M_i}}{(y_i/x_j)_{M_i}}
\notag \\
&\quad\times\frac{\prod_{i,j=1}^n(qx_ix_j)_{N_i}
\prod_{1\leq i<j\leq m}(y_iy_j)_{M_i+M_j}}
{\prod_{i,j=1}^m(y_iy_j)_{M_i}
\prod_{1\leq i<j\leq n}(qx_ix_j)_{N_i+N_j}}
\frac{\prod_{i=1}^m(by_i,cy_i,dy_i,q^{1-M_i}/y_ie)_{M_i}}
{\prod_{i=1}^n(qx_i/b,qx_i/c,qx_i/d,q^{-N_i}e/x_i)_{N_i}}
\notag \\
&\quad\times\sum_{k_1,\dots,k_m=0}^{M_1,\dots,M_m}
\Bigg(\frac{\Delta^{\mathrm{C}}(q^{-1/2}yq^k)}
{\Delta^{\mathrm{C}}(q^{-1/2}y)}
\prod_{i=1}^m\frac{(y_i/b,y_i/c,y_i/d,y_i/e)_{k_i}\, q^{k_i}}
{(by_i,cy_i,dy_i,ey_i)_{k_i}} \notag \\
&\qquad\qquad\quad\times
\prod_{i=1}^m \prod_{j=1}^n
\frac{(q^{N_j}y_ix_j,y_i/x_j)_{k_i}}
{(y_ix_j,q^{-N_j}y_i/x_j)_{k_i}} 
\prod_{i,j=1}^m\frac{(q^{-M_j}y_i/y_j,q^{-1}y_iy_j)_{k_i}}
{(qy_i/y_j,q^{M_j}y_iy_j)_{k_i}}\Bigg).
\label{cst1}
\end{align}

Recalling the notation \eqref{ds}, we have the following discrete 
analogue of \eqref{cit2}:
\begin{multline}\label{wt}
\V{12}{11}{n}\big(a;b,c,d,e,f,g,q^{-N}\big)\\[1mm]
=\frac{(aq,aq/ef,\lambda q/e,\lambda q/f)_{(N^n)}}
{(\lambda q,\lambda q/ef,a q/e,aq/f)_{(N^n)}} \,
\V{12}{11}{n}\bigg(\lambda;\frac{\lambda b}{a},\frac{\lambda c}{a},
\frac{\lambda d}{a},e,f,g,q^{-N}\bigg),
\end{multline}
where $bcdefgt^{n-1}=a^3q^{N+2}$ and $\lambda=a^2q/bcd$.

Finally, the following Karlsson--Minton-type transformation
is an analogue of \eqref{cit3}: 
\begin{align}
\V{2r+8}{2r+7}{n}&\bigg(a;bt^{1-n},\frac{aq^{-M}}{b},c_1q^{k_1},
\dots,c_rq^{k_r},\frac{aq}{c_1},\dots,\frac{aq}{c_r},q^{-N}\bigg)
\notag \\
&=\frac{(aq,t^{n-1}q)_{(N^n)}}
{(bq,t^{n-1}aq/b)_{(N^n)}}\,
\frac{(bq,t^{n-1}bq/a)_{(M^n)}}{(b^2q/a,t^{n-1}q)_{(M^n)}}
\prod_{i=1}^r\frac{(bc_i/a,t^{n-1}c_i/b)_{(k_i^n)}}
{(c_i,t^{n-1}c_i/a)_{(k_i^n)}} \notag \\
&\qquad \times\,
\V{2r+8}{2r+7}{n}\bigg(\frac{b^2}a;bt^{1-n},\frac{bq^{-N}}a,
\frac{bc_1q^{k_1}}{a},\dots,\frac{bc_rq^{k_r}}{a},
\frac{bq}{c_1},\dots,\frac{bq}{c_r},q^{-M}\bigg),
\label{bt}
\end{align}
where the $k_i$ are non-negative integers such that $k_1+\dots+k_r=M+N$.

\subsection{Notes}

For $p=0$ the $\mathrm{A}_n$ summations \eqref{as1}, \eqref{as2} and 
\eqref{as4} are due to Milne \cite{milne97}, Schlosser \cite{schlosser97} 
(see also also \cite{bhatnagar99}), and Gustafson and Rakha
\cite{gustrakha00}, respectively.
The $p=0$ case of the $\mathrm{C}_n$ summation \eqref{dgs} was found
independently by Denis and Gustafson \cite{denisgust92}, and Milne 
and Lilly \cite{milnelilly95}.
The $p=0$ case of the $\mathrm{C}_n$ summation \eqref{ss} is due to 
Schlosser \cite{schlosser00}.
The $p=0$ case of the transformation \eqref{ast1} was obtained, again
independently, by Denis and Gustafson \cite{denisgust92}, and Milne and
Newcomb \cite{milnenewcomb96}. 
The  $p=0$ cases of \eqref{acst} and \eqref{te} are due to Bhatnagar
and Schlosser \cite{bhatschloss98} and Kajihara \cite{kaji04}, 
respectively.
The $p=0$ instances of \eqref{w}, \eqref{bs}, \eqref{cst1}, \eqref{wt} 
and \eqref{bt} were not known prior to the elliptic cases.

For general $p$, the $\mathrm{A}_n$ summations \eqref{as1} and \eqref{as2}
were first obtained by Rosengren \cite{ros04} using an elementary 
inductive argument.
A derivation of \eqref{as1} from \eqref{ai1} using residue calculus is 
given in \cite{spir04} and a similar derivation of \eqref{as2} from 
\eqref{ai3} in \cite{spirwar06}. 
The summation \eqref{as4} was conjectured by Spiridonov \cite{spir04} and 
proved, independently, by Ito and Noumi \cite{itonoumi17c} and 
by Rosengren \cite{ros17}.

As mentioned in the introduction, Warnaar \cite{warnaar02} conjectured
the $\mathrm{C}_n$ summation \eqref{w}.
He also proved the more elementary $\mathrm{C}_n$ summation \eqref{ss}.
Van Diejen and Spiridonov \cite{vandspir00,vandspir01b} showed that the 
$\mathrm{C}_n$ summations \eqref{dgs} and \eqref{w} follow from the
(at that time conjectural) integral identities \eqref{ci1} and \eqref{ci2}. 
This in particular implied the first proof of \eqref{w} for $p=0$.
For general $p$, the summations \eqref{dgs} and \eqref{w} were proved by
Rosengren \cite{ros01,ros04}, using the case $N=1$ of Warnaar's identity
\eqref{ss}. 
Subsequent proofs of \eqref{w} were given in 
\cite{coskgust06,itonoumi17,rains06,rains10}.
The proofs in \cite{coskgust06,rains06} establish the more general
sum \eqref{binomialconvolution} for elliptic binomial coefficients.
In \cite{rains10} the identity \eqref{w} arises as a special case
of the discrete biorthogonality relation \eqref{discretebiorthogonality}
for the elliptic biorthogonal functions $\tilde{R}_{\lambda}$.

The transformations \eqref{ast1} and \eqref{acst} were obtained by
Rosengren \cite{ros04}, together with two more $\mathrm{A}_n$
transformations that are not surveyed here. 
The transformation \eqref{te} was obtained independently by Kajihara
and Noumi \cite{kaji03} and Rosengren \cite{ros06}.
Both these papers contain further transformations that can be obtained
by iterating \eqref{te}. 
The transformation \eqref{cst1} was proved by Rains (personal
communication, 2003) by specializing the parameters of
\cite[Theorem 7.9]{rains10} to a union of geometric progressions. 
It appeared explicitly in \cite[Theorem 4.2]{kmn16} using a similar 
approach to Rains.
The transformation \eqref{wt} was conjectured by Warnaar \cite{warnaar02}
and established by Rains \cite{rains06} using the symmetry of the expression
\eqref{Ericsymmetry} below.
The transformation \eqref{bt} is stated somewhat implicitly by Rains
\cite{rains12}; it includes \eqref{bs} as a special case.

A discrete analogue of the type II $\mathrm{A}_n$ beta integral \eqref{ai2}
has been conjectured by Spiridonov and Warnaar in \cite{spirwar11}.
Surprisingly, this conjecture contains the $\mathrm{C}_n$ identity 
\eqref{w} as a special case.

The summation formula \eqref{ss} can be obtained as a determinant of 
one-dimensional summations. Further summations and transformations of
determinantal type are given in \cite{rosschloss03}.
The special case $t=q$ of \eqref{w} and \eqref{wt} is also closely 
related to determinants, see \cite{schlosser07}.

Transformations related to the sum \eqref{as4} are discussed in \cite{ros17}.
In their work on elliptic Bailey lemmas on root systems, Bhatnagar and 
Schlosser \cite{bhatschloss18} discovered two further elliptic Jackson 
summations for $\mathrm{A}_n$, as well as corresponding transformation
formulas. For none of these an integral analogue is known.
Langer, Schlosser and Warnaar \cite{lsw09} proved a curious $\mathrm{A}_n$
transformation formula, which is new even in the one-variable case. 

\section{Elliptic Macdonald--Koornwinder theory}\label{Sec_EMK}

A function $f$ on $(\mathbb{C}^{\ast})^n$ is said to be
$\mathrm{BC}_n$-symmetric if it is invariant under the action of the
hyperoctahedral group $(\mathbb{Z}/2\mathbb{Z})\wr S_{\hspace{-1pt}n}$.
Here the symmetric group $S_{\hspace{-1pt}n}$ acts by permuting the
variables and $\mathbb{Z}/2\mathbb{Z}$ by replacing a variable with its
reciprocal.
The interpolation functions 
\begin{equation}\label{R-def}
R_{\lambda}^{\ast}(x_1,\dots,x_n;a,b;q,t;p),
\end{equation}
introduced independently by Rains \cite{rains06,rains10} and by Coskun 
and Gustafson \cite{coskgust06}, are $\mathrm{BC}_n$-symmetric elliptic 
functions that generalize Okounkov's $\mathrm{BC}_n$ interpolation
Macdonald polynomials \cite{okounkov98} as well as the Macdonald 
polynomials of type $\mathrm{A}$ \cite{macd95}. 
They form the building blocks of Rains' more general 
$\mathrm{BC}_n$-symmetric functions \cite{rains06,rains10}
\[
\tilde{R}_{\lambda}(x_1,\dots,x_n;a:b,c,d;u,v;q,t;p).
\]
The $\tilde{R}_{\lambda}$ are an elliptic generalization of the Koornwinder 
polynomials \cite{koornwinder92}, themselves a generalization to 
$\mathrm{BC}_n$ of the Askey--Wilson polynomials \cite{askwil85}.
The price one pays for ellipticity is that the functions
$R^{\ast}_{\lambda}$ and $\tilde{R}_{\lambda}$
are neither polynomial nor orthogonal. 
The latter do however form a biorthogonal family, and for $n=1$
they reduce to the continuous biorthogonal functions of Spiridonov 
(elliptic case) \cite{spir04} and Rahman (the $p=0$ case) \cite{rahman86}
and, appropriately specialized, to the discrete biorthogonal functions 
of Spiridonov and Zhedanov (elliptic case) \cite{sz00} and Wilson (the
$p=0$ case) \cite{wilson91}.

\medskip

There are a number of ways to define the elliptic interpolation functions.
Here we will describe them via a branching rule.
The branching coefficient $c_{\lambda\mu}$ is a complex function on
$(\mathbb{C}^{\ast})^7$, indexed by a pair of partitions $\lambda,\mu$.
It is defined to be zero unless $\lambda\succ\mu$, in which case
\begin{align}
c_{\lambda\mu}(z;a,b;q,t,T;p)&=
\frac{(aTz^{\pm},pqa/bt)_{\lambda}}{(aTz^{\pm},pqa/bt)_{\mu}}\,
\frac{(pqz^{\pm}/bt,T)_{\mu}}{(pqz^{\pm}/b,tT)_{\lambda}} \notag \\ 
&\quad \times
\prod_{\substack{(i,j)\in \lambda \\ \lambda'_j=\mu'_j}}
\frac{\theta(q^{\lambda_i+j-1} t^{2-i-\lambda'_j} aT/b)}
     {\theta(p q^{\mu_i-j+1} t^{\mu'_j-i})}
\prod_{\substack{(i,j)\in \lambda \\ \lambda'_j\ne \mu'_j}}
\frac{\theta(q^{\lambda_i-j} t^{\lambda'_j-i+1})}
     {\theta(pq^{\mu_i+j} t^{-i-\mu'_j} aT/b)} \notag \\
&\quad\times \prod_{\substack{(i,j)\in \mu \\ \lambda'_j=\mu'_j}}
\frac{\theta(pq^{\lambda_i-j+1}t^{\lambda'_j-i})}
     {\theta(q^{\mu_i+j-1}t^{1-i-\mu'_j} aT/b)}
\prod_{\substack{(i,j)\in \mu \\ \lambda'_j\ne \mu'_j}}
\frac{\theta(pq^{\lambda_i+j} t^{1-i-\lambda'_j} aT/b)}
     {\theta(q^{\mu_i-j} t^{\mu'_j-i+1})}.
\label{branching}
\end{align}
From \eqref{quasiperiodicity} and the invariance under the substitution
$z\mapsto z^{-1}$ it follows that $c_{\lambda\mu}$ is a 
$\mathrm{BC}_1$-symmetric elliptic function of $z$.
The elliptic interpolation functions are uniquely determined by the 
branching rule
\begin{equation}\label{branchingrule}
R^{*}_{\lambda}(x_1,\dots,x_{n+1};a,b;q,t;p)
=\sum_{\mu} c_{\lambda\mu}(x_{n+1};a,b;q,t,t^n;p) 
R^{*}_{\mu}(x_1,\dots,x_n;a,b;q,t;p),
\end{equation}
subject to the initial condition
$R^{*}_{\lambda}(\text{--}\,;a,b;q,t;p)=
\delta_{\lambda,0}$.
It immediately follows that the interpolation function \eqref{R-def}
vanishes if $l(\lambda)>n$.
From the symmetry and ellipticity of the branching coefficient it also
follows that the interpolation functions are $\mathrm{BC}_1$-symmetric
and elliptic in each of the $x_i$. $S_{\hspace{-1pt}n}$-symmetry (and
thus $\mathrm{BC}_n$-symmetry), however, is not manifest and is a
consequence of  the non-trivial fact that
\begin{equation}\label{zwsymmetry}
\sum_{\mu} c_{\lambda\mu}(z;a,b;q,t,T;p) c_{\mu\nu}(w;a,b;q,t,T/t;p)
\end{equation}
is a symmetric function in $z$ and $w$; see also the discussion around
\eqref{iterated} below.

In the remainder of this section $x=(x_1,\dots,x_n)$.
Comparison of their respective branching rules shows that Okounkov's
$\mathrm{BC}_n$ interpolation Macdonald polynomials
$P^{\ast}_{\hspace{-1pt}\lambda}(x;q,t,s)$ and the ordinary Macdonald
polynomials $P_{\hspace{-1pt}\lambda}(x;q,t)$ arise in the limit as
\[
P_{\hspace{-1pt}\lambda}^{\ast}(x;q,t,s)=
\lim_{p\to 0} \big({-}s^2t^{2n-2}\big)^{-\abs{\lambda}}
q^{-n(\lambda')} t^{2n(\lambda)}
\frac{(t^n)_{\lambda}}{C^{-}_{\lambda}(t)}\, 
R_{\lambda}^{\ast} \big(st^{\delta}x;s,p^{1/2}b;q,t;p\big),
\]
and 
\[
P_{\lambda}(x;q,t)
=\lim_{z\to\infty} z^{-\abs{\lambda}} \lim_{p\to 0} 
\big({-}at^{n-1}\big)^{-\abs{\lambda}} q^{-n(\lambda')} t^{2n(\lambda)}
\frac{(t^n)_{\lambda}}{C^{-}_{\lambda}(t)} \,
R_{\lambda}^{\ast}(zx;a,p^{1/2}b;q,t;p),
\]
where $\delta=(n-1,\dots,1,0)$ is the staircase partition of length
$n-1$, $st^{\delta}x=(st^{n-1}x_1,\dots,st^0 x_n)$ and
$zx=(zx_1,\dots,zx_n)$.

Many standard properties of $P_{\hspace{-1pt}\lambda}(x;q,t)$ and
$P_{\hspace{-1pt}\lambda}^{\ast}(x;q,t,s)$ have counterparts for the
elliptic interpolation functions.
Here we have space for only a small selection.
Up to normalization, Okounkov's $\mathrm{BC}_n$ interpolation Macdonald
polynomials are uniquely determined by symmetry and vanishing properties. 
The latter carry over to the elliptic case as follows:
\begin{equation}\label{vanishing}
R_{\mu}^{\ast}(aq^{\lambda}t^{\delta};a,b;q,t;p)=0
\end{equation}
if $\mu\not\subset\lambda$.
For $q,t,a,b,c,d\in\mathbb{C}^{*}$ the elliptic difference operator 
$D^{(n)}(a,b,c,d;q,t;p)$, acting on $\mathrm{BC}_n$-symmetric functions,
is given by
\[
\big(D^{(n)}(a,b,c,d;q,t;p)f\big)(x)
=\sum_{\sigma\in\{\pm 1\}^n} 
f\big(q^{\sigma/2}x\big)
\prod_{i=1}^n 
\frac{\theta\big(a x_i^{\sigma_i},b x_i^{\sigma_i},
c x_i^{\sigma_i},d x_i^{\sigma_i})}{\theta\big(x_i^{2\sigma_i}\big)}
\prod_{1\leq i<j\leq n} 
\frac{\theta\big(t x_i^{\sigma_i} x_j^{\sigma_j}\big)}
{\theta\big(x_i^{\sigma_i} x_j^{\sigma_j}\big)},
\]
where $q^{\sigma/2}x=(q^{\sigma_1/2}x_1,\dots,q^{\sigma_n/2}x_n)$.
Then
\begin{multline}\label{differenceqn}
D^{(n)}(a,b,c,d;q,t;p) 
R_{\lambda}^{\ast}\big(x;aq^{1/2},bq^{1/2};q,t;p\big) \\
=\prod_{i=1}^n \theta\big(abt^{n-i},acq^{\lambda_i}t^{n-i},
bcq^{-\lambda_i}t^{i-1}\big)
\cdot R_{\lambda}^{\ast}(x;a,b;q,t;p)
\end{multline}
provided that $t^{n-1}abcd=p$.
Like the Macdonald polynomials, there is no simple closed-form expression
for the elliptic interpolation functions. 
When indexed by rectangular partitions of length $n$, however, they do
admit a simple form, viz.\
\begin{equation}\label{rectangle}
R^{*}_{(N^n)}(x;a,b;q,t;p)
=\prod_{i=1}^n \frac{(ax_i^{\pm})_N}{(pqx_i^{\pm}/b)_N}.
\end{equation}
The principal specialization formula for the elliptic interpolation
functions is
\begin{equation}\label{Pspec}
R^{*}_{\lambda}(zt^{\delta};a,b;q,t;p)
=\frac{(t^{n-1}az,a/z)_{\lambda}}{(pqt^{n-1}z/b,pq/bz)_{\lambda}}.
\end{equation}
The $R^{*}_{\lambda}$ satisfy numerous symmetries, all direct consequence
of symmetries of the branching coefficients $c_{\lambda\mu}$.
Two of the most notable ones are
\begin{subequations}
\begin{align}\label{symminus}
R^{*}_{\lambda}(x;a,b;q,t;p)
&=R^{*}_{\lambda}(-x;-a,-b;q,t;p) \\
&=\bigg(\frac{qt^{n-1}a}{b}\bigg)^{2\abs{\lambda}}
q^{4n(\lambda')} t^{-4n(\lambda)} 
R^{*}_{\lambda}(x;1/a,1/b;1/q,1/t;p).
\label{symreciprocal}
\end{align}
\end{subequations}

Specializations of $R^{*}_{\mu}$ give rise to elliptic binomial
coefficients.
Before defining these we introduce the function
\[
\Delta_{\lambda}(a\vert b_1,\dots,b_k)=
\frac{(pqa)_{2\lambda^2}}{C^{-}_{\lambda}(t,pq)C^{+}_{\lambda}(a,pqa/t)}\,
\frac{(b_1,\dots,b_k)_{\lambda}}{(pqa/b_1,\dots,pqa/b_k)_{\lambda}},
\]
where the dependence on $q,t$ and $p$ has been suppressed and 
where $2\lambda^2$ is shorthand for the partition 
$(2\lambda_1,2\lambda_1,2\lambda_2,2\lambda_2,\dots)$.
Explicitly, for $\lambda$ such that $l(\lambda)\leq n$,
\begin{multline*}
\Delta_{\lambda}(a\vert b_1,\dots,b_k)=
\bigg(\frac{(-1)^k a^{k-3} q^{k-3}t}{b_1\cdots b_k}\bigg)^{\abs{\lambda}}
q^{(k-4)n(\lambda')} t^{-(k-6)n(\lambda)}\,
\frac{(at^{1-n},aqt^{-n},b_1,\dots,b_k)_{\lambda}}
{(qt^{n-1},t^n,aq/b_1,\dots,aq/b_k)_{\lambda}}\\
\times
\prod_{i=1}^n \frac{\theta(at^{2-2i}q^{2\lambda_i})}{\theta(at^{2-2i})}
\prod_{1\leq i<j\leq n} \bigg(
\frac{\theta(t^{j-i}q^{\lambda_i-\lambda_j},
at^{2-i-j}q^{\lambda_i+\lambda_j})}
{\theta(t^{j-i},at^{2-i-j})}\,
\frac{(t^{j-i+1})_{\lambda_i-\lambda_j}(at^{3-i-j})_{\lambda_i+\lambda_j}}
{(qt^{j-i-1})_{\lambda_i-\lambda_j}(aqt^{1-i-j})_{\lambda_i+\lambda_j}}\bigg),
\end{multline*}
so that
\begin{equation}\label{V-Delta}
\V{r+1}{r}{n}(a;b_1,\dots,b_{r-4})=\sum_{\lambda} 
\frac{(b_3,\dots,b_{r-4},qt^{n-1}b_1b_2)_{\lambda}}
{(aq/b_3,\dots,aq/b_{r-4},at^{1-n}/b_1b_2)_{\lambda}}\,
\Delta_{\lambda}\bigg(a\Big\vert t^n,b_1,b_2,\frac{at^{1-n}}{b_1b_2}\bigg).
\end{equation}

The elliptic binomial coefficients 
$\binom{\lambda}{\mu}_{[a,b]}=\binom{\lambda}{\mu}_{[a,b];q,t;p}$
may now be defined as
\begin{equation}\label{ellipticbinomial} 
\binom{\lambda}{\mu}_{[a,b]}=
\Delta_{\mu}(a/b\vert t^n,1/b)\,
R^{*}_{\mu}\big(x_1,\dots,x_n;a^{1/2}t^{1-n},ba^{-1/2};q,t;p\big)
\big|_{x_i=a^{1/2}q^{\lambda_i} t^{1-i}}\, ,
\end{equation}
where on the right $n$ can be chosen arbitrarily provided that 
$n\geq l(\lambda),l(\mu)$. 
Apart from their $n$-independence, the elliptic binomial coefficients
are also independent of the choice of square root of $a$.
Although $\binom{\lambda}{0}_{[a,b]}=1$ they are not normalized like
ordinary binomial coefficients, and
\begin{equation}\label{lambdalmabda}
\binom{\lambda}{\lambda}_{[a,b]}=
\frac{(1/b,pqa/b)_{\lambda}}{(b,pqa)_{\lambda}}\,
\frac{C^{+}_{\lambda}(a)}{C^{+}_{\lambda}(a/b)}.
\end{equation}
The elliptic binomial coefficients vanish unless $\mu\subset\lambda$,
are elliptic in both $a$ and $b$, and invariant under the simultaneous
substitution $(a,b,q,t)\mapsto (1/a,1/b,1/q,1/t)$.
They are also conjugation symmetric:
\begin{equation}\label{conjugationsymmetry}
\binom{\lambda}{\mu}_{[a,b];q,t;p}=
\binom{\lambda'}{\mu'}_{[aq/t,b];1/t,1/q;p}.
\end{equation}
A key identity is
\begin{equation}\label{binomialconvolution}
\binom{\lambda}{\nu}_{[a,c]}=
\frac{(b,ce,cd,bde)_{\lambda}}{(cde,bd,be,c)_{\lambda}}\,
\frac{(1/c,bd,be,cde)_{\nu}}{(bcde,e,d,b/c)_{\nu}}
\sum_{\mu} \frac{(c/b,d,e,bcde)_{\mu}}{(bde,ce,cd,1/b)_{\mu}}\,
\binom{\lambda}{\mu}_{[a,b]} \binom{\mu}{\nu}_{[a/b,c/b]}
\end{equation}
for generic parameters such that $bcde=aq$.
The $c\to 1$ limit of
$\binom{\lambda}{\nu}_{[a,c]} (c)_{\lambda}/(1/c)_{\nu}$ exists and is
given by $\delta_{\lambda\nu}$.
Multiplying both sides of \eqref{binomialconvolution} by
$(c)_{\lambda}/(1/c)_{\nu}$ and then letting $c$ tend to $1$ thus yields
the orthogonality relation
\begin{equation}\label{orthogonality}
\sum_{\mu} \binom{\lambda}{\mu}_{[a,b]}\binom{\mu}{\nu}_{[a/b,1/b]}
=\delta_{\lambda\mu}.
\end{equation}

For another important application of \eqref{binomialconvolution} 
we note that by \eqref{Pspec} 
\[
\binom{(N^n)}{\mu}_{[a,b]}
=\Delta_{\mu}\big(a/b\vert t^n,1/b,aq^Nt^{1-n},q^{-N}\big).
\]
Setting $\lambda=(N^n)$ and $\nu=0$ in \eqref{binomialconvolution}
and recalling \eqref{V-Delta} yields the $\mathrm{C}_n$ Jackson 
summation \eqref{w}.
Also \eqref{ss} may be obtained as a special case of
\eqref{binomialconvolution} but the details of the derivation are more
intricate, see \cite{rains06}.
As a final application of \eqref{binomialconvolution} it can be shown
that
\begin{equation}\label{Ericsymmetry}
\frac{(b,b'e)_{\lambda}}{(b'de,bd)_{\lambda}}\,
\frac{(b'de,bf/c)_{\nu}}{(b/c,b'de/f)_{\nu}}
\sum_{\mu} \frac{(c/b,d,e,bb'de,b'g/c,b'f/c)_{\mu}}
{(bb'de/c,b'e,b'd,1/b,f,g)_{\mu}}\,
\binom{\lambda}{\mu}_{[a,b]}
\binom{\mu}{\nu}_{[a/b,c/b]}
\end{equation}
is symmetric in $b$ and $b'$, where $bb'de=aq$ and $cde=fg$.
Setting $\lambda=(N^n)$ and $\nu=0$ results in the transformation formula 
\eqref{wt}.
Now assume that $bcd=b'c'd'$.
Twice using the symmetry of \eqref{Ericsymmetry} it follows that 
\begin{multline}\label{iterated}
\sum_{\mu} \frac{(c,d,aq/c',aq/d')_{\lambda}}{(c,d,aq/c',aq/d')_{\mu}}\,
\frac{(c'/b,d'/b,aq/bc,aq/bd)_{\mu}}{(c'/b,d'/b,aq/bc,aq/bd)_{\nu}}\\
\times
\frac{(1/b,aq/b)_{\nu}}{(1/b,aq/bb')_{\mu}}\,
\frac{(b',aq)_{\mu}}{(b',aq/b)_{\lambda}}\,
\binom{\lambda}{\mu}_{[a,b]}\binom{\mu}{\nu}_{[a/b,b']}
\end{multline}
is invariant under the simultaneous substitution
$(b,c,d)\leftrightarrow (b',c',d')$.
The branching coefficient \eqref{branching} may be expressed as an
elliptic binomial coefficient as
\begin{equation}\label{branchingbinomial}
c_{\lambda\mu}(z;a,b;q,t,T;p)=
\frac{(aTz^{\pm},pqa/bt,t)_{\lambda}}{(aTz^{\pm},pqa/bt,1/t)_{\mu}}\,
\frac{(pqz^{\pm}/bt,T,pqaT/b)_{\mu}}{(pqz^{\pm}/b,tT,pqaT/bt)_{\lambda}}\,
\binom{\lambda}{\mu}_{[aT/b,t]},
\end{equation}
so that, up to a simple change of variables and the use of 
\eqref{quasiperiodicity}, the $z,w$-symmetry of \eqref{zwsymmetry}
corresponds to the $b=b'$ case of the symmetry of \eqref{iterated}.
To conclude our discussion of the elliptic binomial coefficients we
remark that they also arise as connection coefficients between
the interpolation functions.
Specifically,
\begin{multline}\label{connection}
R^{\ast}_{\lambda}(x;a,b;q,t;p) \\
=\sum_{\mu} \binom{\lambda}{\mu}_{[t^{n-1}a/b,a/a']}\,
\frac{(a/a',t^{n-1}aa')_{\lambda}}{(a'/a,t^{n-1}aa')_{\mu}}\,
\frac{(pqt^{n-1}a/b,pq/ab)_{\mu}}{(pqt^{n-1}a'/b,pq/a'b)_{\lambda}}\,
R^{\ast}_{\mu}(x;a',b;q,t;p).
\end{multline}

\medskip

Let $a,b,c,d,u,v,q,t$ be complex parameters such that
$t^{2n-2}abcduv=pq$, and $\lambda$ a partition of length at most $n$. 
Then the $\mathrm{BC}_n$-symmetric biorthogonal functions 
$\tilde{R}_{\lambda}$ are defined as
\begin{multline}\label{Rtildedef}
\tilde{R}_{\lambda}(x;a\hspace{1pt}{:}\hspace{1pt}b,c,d;u,v;q,t;p) \\
=\sum_{\mu\subset\lambda}
\binom{\lambda}{\mu}_{[1/uv,t^{1-n}/av]}\,
\frac{(pq/bu,pq/cu,pq/du,pq/uv)_{\mu}}
{(t^{n-1}ab,t^{n-1}ac,t^{n-1}ad,t^{n-1}av)_{\mu}}\,
R^{*}_{\mu}(x;a,u;q,t;p).
\end{multline}
By \eqref{orthogonality} this relation between the two families of
$\mathrm{BC}_n$ elliptic functions can be inverted.
We also note that from \eqref{connection} it follows that
\begin{equation}\label{reduction}
\tilde{R}_{\lambda}
\big(x;a\hspace{1pt}{:}\hspace{1pt}b,c,d;u,t^{1-n}/b;q,t;p\big)
=\frac{(pq/au,pqt^{n-1}a/u)_{\lambda}}{(b/a,t^{n-1}ab)_{\lambda}}\,
R^{*}_{\mu}(x;b,u;q,t;p),
\end{equation}
so that the interpolation functions are a special case of the biorthogonal
functions.
Finally, from Okounkov's binomial formula for Koornwinder polynomials
\cite{okounkov98} it follows that in the $p\to 0$ limit the
$\tilde{R}_{\lambda}$ simplify to the Koornwinder polynomials
$K_{\lambda}(x;a,b,c,d;q,t)$:
\begin{multline*}
K_{\lambda}\big(x;a,b,c,d;q,t)=
\lim_{p\to 0} \big(at^{n-1}\big)^{-\abs{\lambda}} t^{n(\lambda)}
\frac{(t^n,t^{n-1}ab,t^{n-1}ac,t^{n-1}ad)_{\lambda}}
{C^{-}_{\lambda}(t)C^{+}_{\lambda}(abcdt^{2n-2}/q)} \\
\times \tilde{R}_{\lambda}\big(x;a\hspace{1pt}{:}\hspace{1pt}b,c,d;
up^{1/2},vp^{1/2},q,t;p\big).
\end{multline*}

Most of the previously-listed properties of the interpolation functions
have implications for the biorthogonal functions.
For example, using \eqref{Pspec} and \eqref{binomialconvolution} one can 
prove the principal specialization formula
\begin{equation}\label{PrincipalSpecialization}
\tilde{R}_{\lambda}(bt^{\delta};a\hspace{1pt}{:}\hspace{1pt}b,c,d;u,v;q,t;p)
=\frac{(t^{n-1}bc,t^{n-1}bd,t^{1-n}/bv,pqt^{n-1}a/u)_{\lambda}}
{(t^{n-1}ac,t^{n-1}ad,t^{1-n}/av,pqt^{n-1}b/u)_{\lambda}}.
\end{equation}
Another result that carries over is the elliptic difference equation
\eqref{differenceqn}. Combined with \eqref{Rtildedef} it yields
\begin{multline}\label{differenceequation}
D^{(n)}\big(a,u,b,pt^{1-n}/uab;q,t;p\big) \tilde{R}_{\lambda}^{\ast}\big(x;
aq^{1/2}\hspace{1pt}{:}\hspace{1pt}bq^{1/2},cq^{-1/2},dq^{-1/2};
uq^{1/2},vq^{-1/2},q,t;p\big) \\
=\prod_{i=1}^n \theta\big(abt^{n-i},aut^{n-i},but^{n-i}\big)
\cdot \tilde{R}_{\lambda}(x;a\hspace{1pt}{:}\hspace{1pt}b,c,d;
u,v,q,t;p).
\end{multline}

The Koornwinder polynomials are symmetric in the parameters $a,b,c,d$.
From \eqref{Rtildedef} it follows that $\tilde{R}_{\lambda}$ is 
symmetric in $b,c,d$ but the choice of normalization breaks the
full $S_{\hspace{-1pt}4}$ symmetry. Instead,
\begin{equation}\label{parametersymmetry}
\tilde{R}_{\lambda}(x;a\hspace{1pt}{:}\hspace{1pt}b,c,d;u,v;q,t;p)=
\tilde{R}_{\lambda}(x;b\hspace{1pt}{:}\hspace{1pt}a,c,d;u,v;q,t;p)
\tilde{R}_{\lambda}(bt^{\delta};a\hspace{1pt}{:}\hspace{1pt}b,c,d;u,v;q,t;p).
\end{equation}
For partitions $\lambda,\mu$ such that $l(\lambda),l(\mu)\leq n$
the biorthogonal functions satisfy evaluation symmetry:
\begin{equation}\label{evaluationsymmetry}
\tilde{R}_{\lambda}(at^{\delta}q^{\mu};
a\hspace{1pt}{:}\hspace{1pt}b,c,d;u,v;q,t;p)=
\tilde{R}_{\mu}(\hat{a}t^{\delta}q^{\lambda};
\hat{a}\hspace{1pt}{:}\hspace{1pt}\hat{b},\hat{c},\hat{d};
\hat{u},\hat{v};q,t;p),
\end{equation}
where
\[
\hat{a}=\sqrt{abcd/pq},\quad
\hat{a}\hat{b}=ab,\quad
\hat{a}\hat{c}=ac,\quad
\hat{a}\hat{d}=ad,\quad
a\hat{u}=\hat{a}u,\quad
a\hat{v}=v\hat{a}.
\]

Given a pair of partitions $\lambda,\mu$ such that $l(\lambda),l(\mu)\leq n$,
define
\[
\tilde{R}_{\lambda\mu}(x;a\hspace{1pt}{:}\hspace{1pt}b,c,d;u,v;t;p,q) \\
=\tilde{R}_{\lambda}(x;a\hspace{1pt}{:}\hspace{1pt}b,c,d;u,v;p,t;q)
\tilde{R}_{\mu}(x;a\hspace{1pt}{:}\hspace{1pt}b,c,d;u,v;q,t;p).
\]
Note that
$\tilde{R}_{\lambda\mu}(x;a\hspace{1pt}{:}\hspace{1pt}b,c,d;u,v;t;p,q)$ is
invariant under the simultaneous substitutions $\lambda\leftrightarrow\mu$
and $p\leftrightarrow q$.
The functions 
$\tilde{R}_{\lambda\mu}(x;a\hspace{1pt}{:}\hspace{1pt}b,c,d;u,v;t;p,q)$
form a biorthogonal family, with continuous biorthogonality relation
\begin{align}
\kappa_n^{\mathrm{C}} & \int_{C_{\lambda\nu,\mu\omega}}
\tilde{R}_{\lambda\mu}(z_1,\dots,z_n;
t_1\hspace{1pt}{:}\hspace{1pt}t_2,t_3,t_4;t_5,t_6;t;p,q) 
\tilde{R}_{\nu\omega}(z_1,\dots,z_n;
t_1\hspace{1pt}{:}\hspace{1pt}t_2,t_3,t_4;t_6,t_5;t;p,q) \notag \\
&\qquad\qquad\qquad\qquad\qquad\times
\prod_{1\leq i<j\leq n} \frac{\Gamma(tz_i^{\pm}z_j^{\pm})}
{\Gamma(z_i^{\pm}z_j^{\pm})}
\prod_{j=1}^n\frac{\prod_{i=1}^6\Gamma(t_iz_j^{\pm})}
{\Gamma(z_j^{\pm 2})}\,
\frac{\dup z_1}{z_1}\cdots\frac{\dup z_n}{z_n} \notag \\
&=\delta_{\lambda\nu}\delta_{\mu\omega} 
\prod_{m=1}^n\bigg(\,\frac{\Gamma(t^m)}{\Gamma(t)}
\prod_{1\leq i<j\leq 6}\Gamma(t^{m-1}t_it_j)\bigg) \notag \\
&\quad\times
\frac{1}
{\Delta_{\lambda}(1/t_5t_6\vert t^n,t^{n-1}t_0t_1,t^{n-1}t_0t_2,t^{n-1}t_0t_3,
t^{1-n}/t_0t_5,t^{1-n}/t_0t_6;p,t;q)} \notag \\
&\quad\times
\frac{1}
{\Delta_{\mu}(1/t_5t_6\vert t^n,t^{n-1}t_0t_1,t^{n-1}t_0t_2,t^{n-1}t_0t_3,
t^{1-n}/t_0t_5,t^{1-n}/t_0t_6;q,t;p)}. 
\label{biorthogonality}
\end{align}
Here, $C_{\lambda\nu,\mu\omega}$ is a deformation of $\mathbb{T}^n$ 
which separates sequences of poles of the integrand tending to zero 
from sequences tending to infinity. The location of these poles 
depends on the choice of partitions, see \cite{rains10} for details. 
Provided $\abs{t}<1$ and $\abs{t_i}<1$ for $1\leq i\leq 6$ we can take
$C_{00,00}=\mathbb{T}^n$ so that for $\lambda=\mu=\nu=\omega=0$ one
recovers the type $\mathrm{C}_n^{(\text{II})}$ integral \eqref{ci2}.
The summation \eqref{w}, which is the discrete analogue of \eqref{ci2},
follows in a similar manner from the discrete biorthogonality relation
\begin{align}
\sum_{\mu\subset (N^n)} &
\Delta_{\mu}\big(t^{2n-2}a^2\vert 
t^n,t^{n-1}ac,t^{n-1}ad,t^{n-1}au,t^{n-1}av,q^{-N}\big) \notag \\
& \qquad \times
\tilde{R}_{\lambda}(aq^{\mu} t^{\delta};
a\hspace{1pt}{:}\hspace{1pt}b,c,d;u,v;t;p,q) 
\tilde{R}_{\nu}(aq^{\mu} t^{\delta};
a\hspace{1pt}{:}\hspace{1pt}b,c,d;v,u;t;p,q) \notag \\
&=\frac{\delta_{\lambda\nu}}
{\Delta_{\lambda}(1/uv\vert t^n,t^{n-1}ab,t^{n-1}ac,t^{n-1}ad,
t^{1-n}/au,t^{1-n}/av)} \notag \\[1mm]
&\qquad \times
\frac{(b/a,pq/uc,pq/ud,pq/uv)_{(N^n)}}
{(pqt^{n-1}a/u,t^{n-1}bc,t^{n-1}bd,t^{n-1}bv)_{(N^n)}},
\label{discretebiorthogonality}
\end{align}
where $t^{2n-2}abcduv=pq$ and $q^Nt^{n-1}ab=1$.
The discrete biorthogonality can also be lifted to the functions 
$\tilde{R}_{\lambda\mu}$ but since the resulting identity factors into
two copies of \eqref{discretebiorthogonality} --- the second copy with
$q$ replaced by $p$ and $N$ by a second discrete parameter $M$ --- this
is no more general than the above.

The final result listed here is a (dual) Cauchy identity which incorporates
the Cauchy identities for the Koornwinder polynomials, $\mathrm{BC}_n$
interpolation Macdonald polynomials and ordinary Macdonald polynomials:
\begin{multline}\label{Cauchy}
\sum_{\lambda\subset (N^n)} 
\Delta_{\lambda} 
\big(q^{1-2N}/uv\vert t^n,q^{-N},q^{1-N}t^{1-n}/av,a/u\big) \\
\times
\tilde{R}_{\lambda}(x;a\hspace{1pt}{:}\hspace{1pt}b,c,d;q^Nu,q^{N-1}v;q,t;p)
\tilde{R}_{\hat{\lambda}}(y;a\hspace{1pt}{:}\hspace{1pt}b,c,d;
t^nu,t^{n-1}v;t,q;p) \\
=\frac{(a/u,pq^{1-N}/au,pq^{1-N}/bu,pq^{1-N}/cu,pq^{1-N}/du,
pq^{2-2N}/uv)_{(N^n)}}
{(t^{n-1}ab,t^{n-1}ac,t^{n-1}ad,q^{N-1}t^{n-1}av)_{(N^n)}} \\
\times
\prod_{i=1}^n\prod_{j=1}^N \theta(x_i^{\pm} y_j)
\prod_{i=1}^n \frac{1}{(ux_i^{\pm})_m}
\prod_{j=1}^N \frac{1}{(p/uy_j,y_j/u;1/t,p)_n},
\end{multline}
where $x=(x_1,\dots,x_n)$, $y=(y_1,\dots,y_N)$,
$\hat{\lambda}=(n-\lambda'_m,\dots,n-\lambda'_1)$
and $abcduvq^{2m-2}t^{2n-2}=p$.

\subsection{Notes}\label{Sec_EMKnotes}

Instead of $R_{\lambda}^{\ast}(x_1,\dots,x_n,a,b;q,t;p)$,
Rains denotes the $\mathrm{BC}_n$-symmetric interpolation functions 
as $R_{\lambda}^{\ast(n)}(x_1,\dots,x_n,a,b;q,t;p)$, see 
\cite{rains06,rains10,rains12}.
An equivalent family of functions is defined by Coskun and Gustafson in 
\cite{coskgust06} (see also \cite{cosk08}). 
They refer to these as well-poised Macdonald 
functions, denoted as $W_{\lambda}(x_1,\dots,x_n;q,p,t,a,b)$.
The precise relation between the two families is given by
\begin{multline*}
W_{\lambda}\big(x_1/a,\dots,x_n/a;q,p,t,a^2,a/b\big)\\
=\bigg(\frac{t^{1-n}b^2}{q^2}\bigg)^{\abs{\lambda}} 
q^{-2n(\lambda')} t^{2n(\lambda)}
\frac{(t^n)_{\lambda}}{C^{-}_{\lambda}(t)}\,
\frac{(qt^{n-2}a/b;q,t^2;p)_{2\lambda}}
{(qa/tb)_{\lambda}\, C^{+}_{\lambda}(qt^{n-2}a/b)}\,
R^{\ast}_{\lambda}(x_1,\dots,x_n;a,b;q,t;p).
\end{multline*}
Similarly, Rains writes
$\tilde{R}_{\lambda}^{(n)}(x_1,\dots,x_n;{a:b},c,d;u,v;q,t;p)$
for the biorthogonal functions
instead of $\tilde{R}_{\lambda}(x_1,\dots,x_n;{a:b},c,d;u,v;q,t;p)$,
see again \cite{rains06,rains10,rains12}.

The branching rule \eqref{branchingrule} is the $k=1$ instance of 
\cite[Eq.~(4.40)]{rains06} or the $\mu=0$ case of 
\cite[Eq.~(2.14)]{coskgust06}.
The vanishing property \eqref{vanishing} is
\cite[Corollary~8.12]{rains10} combined with \eqref{reduction}, or
\cite[Theorem~2.6]{coskgust06}.
The elliptic difference equation \eqref{differenceqn} is
\cite[Eq.~(3.34)]{rains06}.
The formula \eqref{rectangle} for the interpolation function
indexed by a rectangular partition of length $n$ is the $\lambda=0$ case
of \cite[Eq.~(3.42)]{rains06} or \cite[Corollary 2.4]{coskgust06}.
The principal specialization formula \eqref{Pspec} is
\cite[Eq.~(3.35)]{rains06}.
The symmetry \eqref{symminus} is 
\cite[Eq.~(3.39)]{rains06} and the symmetry \eqref{symreciprocal}
is \cite[Eq.~(3.38)]{rains06} or \cite[Proposition 2.8]{coskgust06}.
The definition of the elliptic binomial coefficients 
\eqref{ellipticbinomial} is due to Rains, see \cite[Eq.~(4.1)]{rains06}.
Coskun and Gustafson define so-called elliptic Jackson coefficients
\[
\omega_{\lambda/\mu}(z;r;a,b)=\omega_{\lambda/\mu}(z;r,q,p;a,b),
\]
see \cite[Eq.~(2.38)]{coskgust06}.
Up to normalization these are the elliptic binomials coefficients:
\[
\omega_{\lambda/\mu}(z;r;a,b)=
\frac{(1/z,az)_{\lambda}}{(qbz,qb/az)_{\lambda}}
\frac{(qbz/r,qb/azr,bq,r)_{\mu}}{(1/z,az,qb/r^2,1/r)_{\mu}}\,
\binom{\lambda}{\mu}_{[b,r]}.
\]
The value of the elliptic binomials \eqref{lambdalmabda} is 
\cite[Eq.~(4.8)]{rains06}. 
It is equivalent to \cite[Eq.~(2.9)]{coskgust06} and also
\cite[Eq.~(4.23)]{coskgust06}.
The conjugation symmetry \eqref{conjugationsymmetry} of the elliptic 
binomial coefficients is \cite[Corollary 4.4]{rains06}.
The summation \eqref{binomialconvolution} is 
\cite[Theorem 4.1]{rains06}, and is equivalent to the ``cocycle identity''
\cite[Eq.~(3.7)]{coskgust06} for the elliptic Jackson coefficients.
The orthogonality relation \eqref{orthogonality} is 
\cite[Corollary 4.3]{rains06} or \cite[Eq.~(4.16)]{coskgust06}.
The symmetry of \eqref{Ericsymmetry} is
\cite[Theorem 4.9]{rains06} or \cite[Eq.~(3.8)]{coskgust06},
and the symmetry of \eqref{iterated} is \cite[Corollary 4.11]{rains06}.
The expression \eqref{branchingbinomial} for the branching coefficients
is a consequence of \cite[Corollary 4.5]{rains06} or
\cite[Lemma 3.11]{coskgust06}.
The connection coefficient identity \eqref{connection} is
\cite[Corollary 4.14]{rains06}. It is equivalent to
the ``Jackson sum'' \cite[Eq.~(3.6)]{coskgust06} for the
well-poised Macdonald functions $W_{\lambda}$.

Definition \eqref{Rtildedef} of the birthogonal functions
is \cite[Eq.~(5.1)]{rains06}, its principal specialization
\eqref{PrincipalSpecialization} is \cite[Eq.~(5.4)]{rains06}
and the difference equation \eqref{differenceequation} is
\cite[Lemma 5.2]{rains06}.
The parameter and evalation symmetries \eqref{parametersymmetry}
and \eqref{evaluationsymmetry} are
\cite[Theorem 5.1]{rains06} and 
\cite[Theorem 5.4]{rains06}, respectively.
The important biorthogonality relation
\eqref{biorthogonality} is a combination of
\cite[Theorem 8.4]{rains10} and \cite[Theorem 8.10]{rains10}.
Its discrete analogue \eqref{discretebiorthogonality} 
is \cite[Theorem 5.8]{rains06}, see also
\cite[Theorem 8.11]{rains06}.
Finally, the Cauchy identity \eqref{Cauchy} is 
\cite[Theorem 5.11]{rains06}.

The $\mathrm{BC}_n$-symmetric interpolation functions satisfy several 
further important identities not covered in the main text,
such as a ``bulk'' branching rule \cite[Theorem 4.16]{rains06}
which extends \eqref{branchingrule}, and 
a generalized Pieri rule \cite[Theorem 4.17]{rains06}. 
In \cite{cosk08} Coskun applies the elliptic binomial coefficients
(elliptic Jackson coefficients in his language)
to formulate an elliptic Bailey lemma of type $\mathrm{BC}_n$.
The interpolation functions further admit a generalization to skew 
interpolation functions \cite{rains12}
\[
\mathcal{R}_{\lambda/\mu}([v_1,\dots,v_{2n}];a,b;q,t;p),
\qquad \mu\subseteq\lambda.
\]
These are elliptic functions, symmetric in the variables $v_1,\dots,v_{2n}$,
such that \cite[Theorem 2.5]{rains12}
\[
R^{*}_{\lambda}(x_1,\dots,x_n;a,b;q,t;p)
=\frac{(pqa/tb)_{\lambda}}{(t^n)_{\lambda}}\,
\mathcal{R}^{\ast}_{\lambda/0}
\big(\big[t^{1/2}x_1^{\pm},\dots,t^{1/2}x_n^{\pm}\big];
t^{n-1/2}a,t^{1/2}b;q,t;p\big).
\]
They also generalize the $n$-variable skew elliptic Jackson coefficients
\cite[Eq.~(2.43)]{coskgust06}
\[
\omega_{\lambda/\mu}(x_1,\dots,x_n;r,q,p;a,b)
\]
of Coskun and Gustafson:
\begin{multline*}
\omega_{\lambda/\mu}\big(r^{n-1/2}x_1/a,\dots,
r^{n-1/2}x_n/a;r;a^2r^{1-2n},ar^{1-n}/b\big) \\
=\bigg({-}\frac{b^3}{q^3a}\bigg)^{\abs{\lambda}-\abs{\mu}} 
q^{3n(\mu')-3n(\lambda')} t^{3n(\lambda)-3n(\mu)} r^{-n\abs{\mu}} 
\frac{(aq/br)_{\lambda}}{(aqr^{-n-1}/b)_{\mu}} \, 
\frac{(r)_{\mu}}{(r)_{\lambda}}  \\ \times
\mathcal{R}^{\ast}_{\lambda/\mu}
\big(\big[r^{1/2}x_1^{\pm},\dots
r^{1/2}x_n^{\pm}\big];a,b;q,t;p\big).
\end{multline*}
A very different generalization of the interpolation functions
is given in \cite{rains18} in the form of an interpolation kernel 
$\mathcal{K}_c(x_1,\dots,x_n;y_1,\dots,y_n;q,t;p)$.
By specialising $y_i=q^{\lambda_i} t^{n-i}a/c$ with $c=\sqrt{t^{n-1}ab}$
for all $1\leq i\leq n$ one recovers, up to a simple normalising factor, 
$R^{\ast}_{\lambda}(x_1,\dots,x_n;a,b;q,t;p)$.
In the same paper Rains uses this kernel to prove quadratic transformation 
formulas for elliptic Selberg integrals.

Also for the biorthogonal functions we have omitted a number of
further results, such as a ``quasi''-Pieri formula 
\cite[Theorem 5.10]{rains06} and a connection coefficient formula of
Askey--Wilson type \cite[Theorem 5.6]{rains06}, 
generalizing \eqref{connection}.

\end{document}